\documentclass[11pt,reqno]{amsart}
\usepackage{amsmath, amsthm, amscd, amssymb, amsfonts, latexsym, graphicx}
\usepackage{fullpage}
\usepackage{bm}
\usepackage[all]{xy}
\usepackage{tikz}
\usepackage{mathrsfs}
\usepackage{xcolor}
\usepackage[hidelinks]{hyperref}

\usetikzlibrary{patterns,decorations.pathreplacing}

\renewcommand{\theequation}{\arabic{section}.\arabic{equation}}
\numberwithin{equation}{section}

\theoremstyle{plain}
\newtheorem{lem}{Lemma}[section]
\newtheorem{theorem}[lem]{Theorem}
\newtheorem{corollary}[lem]{Corollary}

\theoremstyle{definition}

\newtheorem{conjecture}[lem]{Conjecture}

\newcommand{\G}{\Gamma}
\newcommand{\eps}{\varepsilon}
\newcommand{\cB}{\mathcal B}

\newcommand{\sym}{\operatorname{sym}}
\newcommand{\vol}{\operatorname{vol}}
\newcommand{\supp}{\operatorname{supp}}

\title{Quantum Variance for Hecke--Maass Cusp Forms on the Vertical Geodesic}
\author{Qingfeng Sun}
	\address{School of Mathematics and Statistics, Shandong University, Weihai\\Weihai, Shandong 264209, China}
	\address{State Key Laboratory of Cryptography and Digital Economy
Security \\ Shandong University, Weihai, 264209, China}
	\email{qfsun@sdu.edu.cn}

    \author{Qizhi Zhang}
	\address{School of Mathematics and Statistics, Shandong University, Weihai\\Weihai, Shandong 264209, China}
    \email{qzzhang@mail.sdu.edu.cn}

\thanks{Q. Sun was partially supported by the National Natural Science Foundation of China
(Grant No. 12471005) and
the Natural Science Foundation of Shandong Province (Grant No. ZR2023MA003).}

\date{}

\begin{document}
\begin{abstract}
Let $X=\mathrm{PSL}_2(\mathbb Z)\backslash\mathbb H$, and let
$\{u_j\}$ be an orthonormal basis of even Hecke--Maass cusp forms on $X$
with Laplace eigenvalues $1/4+t_j^2$. For a real-valued
$\psi\in C_c^\infty(0,\infty)$, define
\[
 \mu_j(\psi)=\int_0^\infty\psi(y)|u_j(iy)|^2\,\frac{dy}{y},
 \qquad
 M(\psi)=\frac6\pi\int_0^\infty\psi(y)\,\frac{dy}{y}.
\]
For $T\ge2$, $1\le U\le T$, and every $\eps>0$, we prove
\[
 \sum_{\substack{j\ {\rm even}\\T<t_j\le T+U}}
 L(1,\sym^2u_j)|\mu_j(\psi)-M(\psi)|^2
 \ll_{\psi,\eps}T^{1+\eps}.
\]
As a consequence, for fixed $0<\theta\le1$, $0<\eta<\theta/2$, and every
$\eps>0$, if $T^\theta\le U\le T$, then all but
\[
 O_{\psi,\theta,\eta,\eps}
 \!\left(T^{1+2\eta+\eps}\right)
\]
of the even forms with $T<t_j\le T+U$ satisfy
$\mu_j(\psi)=M(\psi)+O(T^{-\eta})$. In particular, the exceptional
proportion tends to zero in this short spectral interval.
\end{abstract}
\maketitle

\section{Introduction}\label{sec:introduction}

The quantum unique ergodicity problem asks whether high-energy Laplace
eigenfunctions become equidistributed on the underlying surface. For the
modular surface $X=\Gamma\backslash\mathbb H$ with  $\mathbb H=\{x+iy:y>0\}$ and
$\Gamma=\mathrm{PSL}_2(\mathbb Z)$, arithmetic QUE is known for
$L^2$-normalised
Hecke--Maass cusp forms by the work of Lindenstrauss
\cite{Lindenstrauss2006} and Soundararajan \cite{Soundararajan2010}.
In particular, for every $\psi\in C_c^\infty(X)$, as $t_j\to\infty$ one has
\[
 \int_X\psi(z)|u_j(z)|^2\,d\mu(z)
 \longrightarrow\frac{1}{\vol(X)}\int_X\psi(z)\,d\mu(z),
\]
where the hyperbolic area measure $d\mu(z)$ is given by
$d\mu(z)=dx\,dy/y^2$, so that $\vol(X)=\pi/3$.

A natural question is whether the restriction mass of these eigenfunctions
also becomes equidistributed along a fixed geodesic. Since the usual QUE
result concerns integration against smooth functions on the surface,
it does not directly
give information about the restriction mass along a geodesic. Young formulated the
following conjecture for the vertical geodesic
\cite[Conjecture~1.1]{Young2016}.

\begin{conjecture}\label{conj:intro-rque}
Let $u_j$ range over even, $L^2$-normalised Hecke--Maass cusp forms on
$X$, with Laplace eigenvalues $1/4+t_j^2$. For every real-valued
$\psi\in C_c^\infty(0,\infty)$,

\begin{equation}\label{eq:intro-rque-prediction}
 \lim_{\substack{t_j\to\infty\\u_j\ \mathrm{even}}}
 \int_0^\infty\psi(y)|u_j(iy)|^2\,\frac{dy}{y}
 =\frac6\pi\int_0^\infty\psi(y)\,\frac{dy}{y}.
\end{equation}
\end{conjecture}

One way to study this conjecture is through quantum variance.
In the holomorphic setting, Zenz \cite{Zenz2021} obtained an asymptotic formula
for the quantum variance along the vertical geodesic, averaged over Hecke
cusp forms of even weights $k\asymp K$. The leading term is of order
$K^{3/2}\log K$. More recently, the authors \cite{SunZhang2026}
extended this result to short weight intervals of length $G=K^{9/10}$,
obtaining an asymptotic formula with a leading term of order
$K^{1/2}G\log K$.

Let
$\{u_j\}$ be an orthonormal basis of even Hecke--Maass cusp forms on $X$
with Laplace eigenvalues $1/4+t_j^2$.
We have the  following Fourier--Whittaker expansion (see \cite{Iwaniec2002})
\[
u_j(x+iy)=\sqrt y\sum_{n\ne0}\rho_j(n)
 K_{it_j}(2\pi|n|y)e(nx),
\]
where, as usual, $e(x)=e^{2\pi i x}$ and
$\rho_j(\pm n)=\rho_j(1)\lambda_j(n)$ for $n\ge1$.
The Hecke eigenvalues $\lambda_j(n)$ satisfy
\[
 \lambda_j(m)\lambda_j(n)
 =\sum_{d\mid(m,n)}\lambda_j(mn/d^2).
\]
The $L^2$-normalisation gives
\[
\frac{|\rho_j(1)|^2}{\cosh(\pi t_j)}
 =\frac{2}{L(1,\sym^2u_j)}.
\]

Our main result is the following.

\begin{theorem}\label{thm:main}
For $T\ge2$, $1\le U\le T$, every real-valued
$\psi\in C_c^\infty(0,\infty)$, and every $\eps>0$,
\begin{equation}\label{eq:mainbound}
 \sum_{\substack{u_j\ \mathrm{even}\\T<t_j\le T+U}}
 L(1,\sym^2u_j)\left|\int_0^\infty\psi(y)|u_j(iy)|^2\,\frac{dy}{y}
 -\frac6\pi\int_0^\infty\psi(y)\,\frac{dy}{y}\right|^2
 \ll_{\psi,\eps}T^{1+\eps}.
\end{equation}
Here the implied constant depends only on $\psi$ and $\eps$, and the estimate is
uniform for $1\le U\le T$.

\end{theorem}

For $U\ge T^\theta$ with fixed $\theta>0$, the even-spectrum
Weyl law (see Lemma~\ref{prop:parity-weyl})
gives $$\#\{j\ \mathrm{even}:T<t_j\le T+U\}\asymp TU.$$
Thus the upper bound in \eqref{eq:mainbound} implies
the average order of the
weighted squares
$L(1,\sym^2u_j)|\mu_j(\psi)-M(\psi)|^2$ is
$O_{\psi,\eps}(T^\eps/U)$. We can compare this bound
with the holomorphic case.
For holomorphic cusp forms of weights $k\asymp K$, the results of
Zenz \cite{Zenz2021} and the authors \cite{SunZhang2026}
give an average of order $K^{-1/2}\log K$, while
for $U=T^\theta$ with $\theta>1/2$, the bound in \eqref{eq:mainbound}
gives
$T^{-\theta+\eps}$, which has a smaller power of the spectral
parameter.

As a direct consequence of Theorem~\ref{thm:main} and the even-spectrum
Weyl law in Lemma~\ref{prop:parity-weyl}, we obtain density-one restricted equidistribution in short
spectral intervals.
\begin{corollary}\label{cor:rque}
Fix a real-valued $\psi\in C_c^\infty(0,\infty)$, $0<\theta\le1$, and
$0<\eta<\theta/2$. If $T^\theta\le U\le T$, then for every $\eps>0$ all
but
\begin{equation}\label{eq:exceptional}
 O_{\psi,\theta,\eta,\eps}(T^{1+2\eta+\eps})
\end{equation}
even forms with $T<t_j\le T+U$ satisfy
\[
 \int_0^\infty\psi(y)|u_j(iy)|^2\,\frac{dy}{y}
 =\frac6\pi\int_0^\infty\psi(y)\,\frac{dy}{y}
 +O(T^{-\eta}).
\]
In particular, taking $0<\eps<\theta-2\eta$, the exceptional proportion
tends to zero.
\end{corollary}

\begin{proof}
Fix $\eps>0$ and let $\mathcal E$ be the exceptional set.
Applying the lower bound
$L(1,\sym^2u_j)\gg_\eps T^{-\eps/2}$ and
Theorem~\ref{thm:main}, we obtain
\[
 \#\mathcal E\ll_{\psi,\eps}T^{1+2\eta+\eps}.
\]
By Lemma~\ref{prop:parity-weyl} in Appendix~\ref{app:parity-weyl},
the even-spectrum count satisfies
\[
 \#\{j\ \mathrm{even}:T<t_j\le T+U\}
 =\frac{2TU+U^2}{24}+O(T\log T)\asymp_\theta TU
\]
uniformly for $T^\theta\le U\le T$ as $T\to\infty$.
Thus by choosing $0<\eps<\theta-2\eta$, we obtain that
the exceptional proportion is
\[
 O_{\psi,\theta,\eta,\eps}(T^{2\eta+\eps-\theta})=o(1).
\]
\end{proof}

\noindent{\bf Outline of the proof of Theorem~\ref{thm:main}.}
We substitute the Fourier--Whittaker expansion of the Maass forms into the
restriction integral and separate the resulting sum into diagonal and
off-diagonal parts.
The diagonal part has a direct
interpretation in terms of an incomplete Eisenstein series. For
$h(y)=y\psi(y)$, let
\[
 E(z\mid h)=\sum_{\gamma\in\Gamma_\infty\backslash\Gamma}
 h(\operatorname{Im}\gamma z).
\]
Unfolding shows that the diagonal part is exactly
\[
 2\langle |u_j|^2,E(\cdot\mid y\psi)\rangle_X.
\]
Since $\vol(X)=\pi/3$, we have
\[
 \frac{2}{\vol(X)}\langle E(\cdot\mid y\psi),1\rangle_X
 =\frac6\pi\int_0^\infty\psi(y)\,\frac{dy}{y}.
\]
Thus the constant part of the incomplete Eisenstein series gives the
predicted mass. The diagonal fluctuations are controlled by the quantum variance
estimate of Luo and Sarnak \cite[Proposition~3.1]{LuoSarnak1995}.

The off-diagonal terms involve the normalised Bessel product integral
\[
 \mathcal B_{m,n}(t)=16\cosh(\pi t)
 \int_0^\infty\psi(y)K_{it}(2\pi my)K_{it}(2\pi ny)\,dy,
 \qquad m\ne n.
\]
In the holomorphic case, the restriction integral involves exponential
functions, which are relatively easy to integrate. For Maass forms, the
corresponding integrals involve $K$-Bessel functions and are harder to
analyse. However, the differential equation for the $K$-Bessel function
allows us to truncate the off-diagonal sum. We also establish estimates for
the spectral derivatives of the corresponding Bessel integrals, which are
needed to apply the spectral large sieve. More precisely, in
Sections~\ref{sec:prelim} and~\ref{sec:trunc}, we show that the
off-diagonal sum can be restricted to
\[
 1\le m<n\ll_\psi T,\qquad n^2-m^2\ll T^{1+\eps},
\]
with a negligible remainder in the off-diagonal contribution.

Applying the Hecke relation and making the change of variables
$m=da$, $n=d(a+b)$, we can rewrite the off-diagonal second moment as
\[
 \sum_{\substack{j\ \mathrm{even}\\T<t_j\le T+U}}
 \frac1{L(1,\sym^2u_j)}\left|
 \sum_{d,a,b\ge1}\lambda_j\bigl(a(a+b)\bigr)
 \mathcal B_{da,d(a+b)}(t_j)\right|^2.
\]
The coefficients in the inner sum still depend on $t_j$.
The spectral
derivative estimates allow us to approximate these Bessel integrals by
finite Fourier series on shorter spectral intervals, with a
negligible error.
Schematically, for an interval $I$ with midpoint $t_I$, this takes the form
\[
 \mathcal B_{m,n}(t)
 \quad\rightsquigarrow\quad
 \sum_{|k|\le K}c_{m,n,I}(k)
 e\!\left(\frac{k(t-t_I)}{|I|}\right),\qquad t\in I.
\]
For each Fourier frequency, the coefficients on the right are independent
of the varying spectral parameter. We may therefore apply Jutila's spectral
large sieve on each short interval. Combining the resulting estimates over
the dyadic ranges proves the off-diagonal bound
$O_{\psi,\eps}(T^{1+\eps})$.

Together with the diagonal variance estimate, this proves
Theorem~\ref{thm:main}.

\medskip


\noindent{\bf Organisation of the paper.}
Section~\ref{sec:prelim} records the spectral normalisations and the
large-sieve input. Section~\ref{sec:trunc} proves the arithmetic truncation,
and Section~\ref{sec:bessel-bounds} develops the Bessel analysis and the
spectral derivative estimates. Section~\ref{sec:sieve} gives the Fourier
localisation needed for the local spectral large sieve. Section~\ref{sec:offdiag}
establishes the dyadic coefficient estimates and completes the proof of the
main theorem. Appendix~\ref{app:parity-weyl} records the parity-specific Weyl
law used in Corollary~\ref{cor:rque}.

\section{Preliminaries}\label{sec:prelim}

Throughout the paper, write
\[
\mu_j(\psi)=\int_0^\infty\psi(y)|u_j(iy)|^2\,\frac{dy}{y},
 \qquad M(\psi)=\frac6\pi\int_0^\infty\psi(y)\,\frac{dy}{y},
 \qquad L_j=L(1,\sym^2u_j),
\]
where $u_j$ is an $L^2$-normalised Hecke--Maass cusp form. We use the
standard bounds of Hoffstein--Lockhart \cite{HoffsteinLockhart1994} and Li
\cite{Li2010},
\[
 t_j^{-\eps}\ll_\eps L_j\ll_\eps t_j^\eps.
\]
Recall the Fourier--Whittaker expansion
\begin{equation}\label{eq:fourier}
u_j(z)=\sum_{n\ne0}\rho_j(n)\sqrt yK_{it_j}(2\pi|n|y)e(nx),
 \qquad
 \rho_j(\pm n)=\rho_j(1)\lambda_j(n)\quad(n\ge1),\end{equation}
where $\|u_j\|_2=1$ gives
\[
\frac{|\rho_j(1)|^2}{\cosh(\pi t_j)}
 =\frac{2}{L(1,\sym^2u_j)}.
\]
The restriction integral $\mu_j(\psi)$ then splits into a diagonal part, treated
using an incomplete Eisenstein series, and an off-diagonal part, reorganised
using the Hecke relations.

For positive integers $m,n$, recall the normalised Bessel product integral
\begin{equation}\label{eq:B-def}
 \cB_{m,n}(t)=16\cosh(\pi t)
 \int_0^\infty\psi(y)K_{it}(2\pi my)K_{it}(2\pi ny)\,dy.
\end{equation}
Substituting the Fourier expansion into the restriction integral gives
\begin{equation}\label{eq:mu-double}
\begin{aligned}
\int_0^\infty\psi(y)|u_j(iy)|^2\,\frac{dy}{y}
&=4|\rho_j(1)|^2
\sum_{m,n\ge1}\lambda_j(m)\lambda_j(n)
\int_0^\infty\psi(y)K_{it_j}(2\pi my)
K_{it_j}(2\pi ny)\,dy \\
&=\frac{1}{2L_j}
\sum_{m,n\ge1}\lambda_j(m)\lambda_j(n)\cB_{m,n}(t_j).
\end{aligned}
\end{equation}
Write $\mu_j=D_j+S_j$, where $D_j$ is the contribution $m=n$.

Let
\begin{equation}\label{eq:incomplete-eis}
 E(z\mid h)=\sum_{\gamma\in\G_\infty\backslash\G}
 h(\Im\gamma z),\qquad h(y)=y\psi(y).
\end{equation}
Unfolding \eqref{eq:incomplete-eis}, using both signs in \eqref{eq:fourier},
gives the exact identity
\begin{equation}\label{eq:diag-eis}
 D_j(\psi)=2\langle |u_j|^2,E(\,\cdot\mid y\psi)\rangle_X.
\end{equation}
Moreover,
\[
 \langle E(\,\cdot\mid y\psi),1\rangle_X
 =\int_0^\infty\psi(y)\,\frac{dy}{y}.
\]
The mean of the incomplete Eisenstein series therefore determines the natural centering constant
\begin{equation}\label{eq:center}
 M(\psi)=\frac{2}{\vol(X)}\int_0^\infty\psi(y)\frac{dy}{y}
 =\frac6\pi\int_0^\infty\psi(y)\frac{dy}{y}.
\end{equation}

\subsection{The diagonal variance estimate}

The diagonal contribution satisfies the following estimate.

\begin{lem}\label{lem:LS-incomplete}
Let $\psi\in C_c^\infty(0,\infty)$ be real-valued.  For $T\ge2$ and
$1\le U\le T$,
\begin{equation}\label{eq:Dvariance}
 \sum_{\substack{j\ {\rm even}\\T<t_j\le T+U}}
 L_j|D_j(\psi)-M(\psi)|^2\ll_{\psi,\eps}T^{1+\eps}.
\end{equation}
\end{lem}

\begin{proof}
By \eqref{eq:diag-eis} and \eqref{eq:center},
\[
 D_j(\psi)-M(\psi)
 =2\left(\langle |u_j|^2,E(\,\cdot\mid y\psi)\rangle_X
 -\frac{\langle E(\,\cdot\mid y\psi),1\rangle_X}{\vol(X)}\right).
\]
Apply \cite[Proposition~3.1]{LuoSarnak1995} to the
incomplete Eisenstein series $E(z\mid h)$ of \eqref{eq:incomplete-eis},
with $h(y)=y\psi(y)$.  We have only replaced the notation $P_{h,0}$
used there by $E(z\mid h)$, and their spectral parameter $\kappa_j$
by $t_j$.  Since $h$ is fixed and compactly supported in $(0,\infty)$,
the derivative constants in that proposition depend only on $\psi$.
It follows that
\[
 \sum_{t_j\le 2T}|D_j(\psi)-M(\psi)|^2
 \ll_{\psi,\eps}T^{1+\eps}.
\]
The asserted bound follows by restricting to the even forms in
$T<t_j\le T+U\le2T$ and using $L_j\ll_\eps T^\eps$,
after adjusting $\eps$ in the variance estimate.
\end{proof}

Using the symmetry in $m,n$ and the Hecke relation
\begin{equation}\label{eq:hecke}
 \lambda_j(m)\lambda_j(n)
 =\sum_{d\mid(m,n)}\lambda_j(mn/d^2)
\end{equation}
and making the substitution $m=da$, $n=d(a+b)$, we obtain
\begin{equation}\label{eq:S-hecke}
 S_j(\psi)=\frac1{L_j}\sum_{d,a,b\ge1}
 \lambda_j(a(a+b))\cB_{da,d(a+b)}(t_j),
\end{equation}

\subsection{Truncation of the off-diagonal sum}

We next make a preliminary truncation of the off-diagonal sum.  Write
$\supp\psi\subset[y_0,y_1]$, where $0<y_0<y_1$, and let
$\mathcal E_j^{\mathrm{far}}$ denote the contribution to $S_j(\psi)$ from
$\max\{m,n\}\gg_\psi T$.

\begin{lem}\label{lem:K-tail}
For every $A>0$, uniformly for $T<t_j\le T+U\le2T$,
\[
 \mathcal E_j^{\mathrm{far}}\ll_{A,\psi}T^{-A}.
\]
Consequently, the off-diagonal sum can be restricted to
\[
 m,n\ll_\psi T.
\]
\end{lem}

\begin{proof}
The integral representation (see \cite[8.432.1]{GradshteynRyzhik2007})
\[
 K_{it}(x)=\int_0^\infty e^{-x\cosh u}\cos(tu)\,du
\]
and the inequality $\cosh u\ge1+u^2/2$ give
\[
 |K_{it}(x)|\le K_0(x)
 \le e^{-x}\int_0^\infty e^{-xu^2/2}\,du
 =\sqrt{\frac{\pi}{2x}}e^{-x}.
\]
Since $y\ge y_0$ on $\supp\psi$ and $t_j\le2T$, it follows that, for
$n\gg_\psi T$,
\[
 e^{\pi t_j/2}|K_{it_j}(2\pi ny)|
 \ll_\psi n^{-1/2}\exp\{\pi T-2\pi y_0n\}
 \ll_{A,\psi}T^{-A}n^{-A}.
\]
The same estimate holds with $m$ in place of $n$.  Substitution in the
off-diagonal part of \eqref{eq:mu-double}, followed by absolute summation,
proves the assertion.
\end{proof}

Fix a sufficiently large constant $C_\psi$ as in Lemma~\ref{lem:K-tail}.
Henceforth we restrict the Fourier indices to $m,n\le C_\psi T$, with a
negligible error.
\subsection{The spectral large sieve}

We record Jutila's spectral large-sieve inequality
\cite[Theorem~1.1]{Jutila2000} in the form needed for the off-diagonal
contribution.

\begin{lem}\label{lem:spectral-sieve}
Let $T\ge2$, $1\le H\le T$, and let $(a_n)$ be a finite complex sequence
supported on $n\le N$.  Then
\begin{equation*}
 \sum_{T\le t_j\le T+H}\frac{1}{L_j}
 \left|\sum_{n\le N}a_n\lambda_j(n)\right|^2
 \ll_\eps (TH+N)(TN)^\eps\sum_{n\le N}|a_n|^2.
\end{equation*}
The sum may be restricted to any subfamily, in particular to the even forms.
\end{lem}

Note that the truncation $m,n\ll_\psi T$ alone is insufficient.
The following calculation explains the need for a further truncation. Fix
$T\le t\le2T$ and consider a dyadic block
\[
 d\asymp D,\qquad a\asymp X,\qquad b\asymp Y,
 \qquad N:=X(X+Y),
\]
for which, uniformly on $\supp\psi$,
\[
 T/4\le 2\pi day\le 2\pi d(a+b)y\le T/2.
\]  Thus both arguments lie in the
oscillatory range at a fixed proportional distance from the turning point.
The standard oscillatory asymptotic for the $K$-Bessel function
\cite[Lemma~3.1]{GhoshReznikovSarnak2013} then gives
\begin{equation*}
 e^{\pi t/2}|K_{it}(2\pi day)|
 +e^{\pi t/2}|K_{it}(2\pi d(a+b)y)|
 \ll_\psi T^{-1/2}.
\end{equation*}
Hence, by
\eqref{eq:B-def},
\begin{equation*}
 \cB_{da,d(a+b)}(t)\ll_\psi T^{-1}.
\end{equation*}

For each $A$, choosing a divisor $a\mid A$ determines $b$ through
$a(a+b)=A$, while there are $O(D)$ choices for $d$.  Thus Cauchy's inequality
and the divisor bound give
\begin{equation*}
 \begin{aligned}
 &\sum_{A\asymp X(X+Y)}
 \left|
 \sum_{\substack{d\asymp D,\ a\asymp X,\ b\asymp Y\\
                  a(a+b)=A}}
 \cB_{da,d(a+b)}(t)
 \right|^2 \\
 &\qquad\ll_{\eps,\psi}
 DT^\eps
 \sum_{\substack{d\asymp D,\ a\asymp X,\ b\asymp Y}}
 |\cB_{da,d(a+b)}(t)|^2
 \ll_{\eps,\psi}
 T^\eps\frac{D^2XY}{T^2}.
 \end{aligned}
\end{equation*}
The spectral large sieve can not be used here directly, but the structure of it would therefore produce the scale
\begin{equation}\label{eq:prelim-large-sieve-scale}
T^\eps\bigl(TU+X(X+Y)\bigr)\frac{D^2XY}{T^2}.
\end{equation}
For $DX\asymp DY\asymp T$, the right-hand side of
\eqref{eq:prelim-large-sieve-scale} is of size
\[
 T^\eps(TU+T^2/D^2).
\]
In particular, blocks with $D\asymp1$ permit the scale
$T^\eps(TU+T^2)$,
which is too large for the desired $T^{1+\eps}$ bound. Thus the
restrictions $m,n\ll_\psi T$ do not suffice. A further restriction on
$d|b|=|m-n|$ is established in Section~\ref{sec:trunc}.

\subsection{Two auxiliary calculus lemmas}

\begin{theorem}[Fa\`a di Bruno formula,  \cite{FaaDiBruno1857}]\label{thm:faa-di-bruno}
Let $f,g\in C^\ell$ and $\ell\ge1$. Then
\begin{equation*}
 \frac{d^\ell}{dh^\ell}f(g(h))
 =\ell!\!\!\sum_{\substack{m_1,\ldots,m_\ell\ge0\\
 m_1+2m_2+\cdots+\ell m_\ell=\ell}}
 f^{(m_1+\cdots+m_\ell)}(g(h))
 \prod_{r=1}^{\ell}\frac{1}{m_r!}
 \left(\frac{g^{(r)}(h)}{r!}\right)^{m_r}.
\end{equation*}
\end{theorem}

\begin{lem}\label{lem:schwartz-majorant}
Let $f:\mathbb R\to[0,\infty)$ be upper semicontinuous and suppose that, for every
$N\ge0$,
\[
 f(x)\ll_N\langle x\rangle^{-N}.
\]
Then there is a nonnegative $\rho\in\mathcal S(\mathbb R)$ such that
$f\le\rho$.
\end{lem}

\begin{proof}
Choose $\chi\in C_c^\infty((-1,1))$ with $\chi\ge0$ and
$\int_{\mathbb R}\chi=1$, and set
\[
 f^\sharp(x):=\sup_{|y-x|\le1}f(y),\qquad
 \rho=f^\sharp*\chi.
\]
Then $f^\sharp(x)\ll_N\langle x\rangle^{-N}$ for every $N$.  Moreover,
if $\chi(x-y)\ne0$, then $|x-y|<1$ and hence $f^\sharp(y)\ge f(x)$.
Therefore
\[
 \rho(x)=\int_{\mathbb R}f^\sharp(y)\chi(x-y)\,dy\ge f(x).
\]
Finally,
\[
 \rho^{(k)}=f^\sharp*\chi^{(k)},
\]
and the compact support of $\chi^{(k)}$, together with the rapid decay of
$f^\sharp$, gives
\[
 |\rho^{(k)}(x)|\ll_{k,N}\langle x\rangle^{-N}
\]
for every $k,N\ge0$.  Thus $\rho\in\mathcal S(\mathbb R)$.
\end{proof}

\section{Truncation of the shift}\label{sec:trunc}

We now prove a further restriction on the off-diagonal sum. Up to a negligible error, the off-diagonal terms may be restricted to
\[
|m^2-n^2|\ll T^{1+\eps}.
\]

The proof uses the fact that the differential operator arising
from the modified Bessel equation is formally symmetric with respect to
the weighted inner product defined below; see, for example,
\cite[\S~5.4, (5.53)]{Teschl2012}.  The modified Bessel function satisfies
\begin{equation}\label{eq:K-bessel-ode}
 x^2K_{it}''(x)+xK_{it}'(x)-(x^2-t^2)K_{it}(x)=0;
\end{equation}
see \cite[8.494(1)]{GradshteynRyzhik2007}.  If $L_x=x\frac{d}{dx}$ denotes the Euler differential operator, then \eqref{eq:K-bessel-ode} becomes
\[
 (L_x^2+t^2)K_{it}(x)=x^2K_{it}(x).
\]
After putting $x=2\pi n y$, we have $x\partial_x=y\partial_y$, and hence
\[
 \bigl((y\partial_y)^2+t^2\bigr)K_{it}(2\pi n y)
   =(2\pi n)^2y^2K_{it}(2\pi n y).
\]
In Sturm--Liouville form, this leads to the differential expression
\[
 P_t:=\frac1y\left(\frac{d}{dy}\,y\frac{d}{dy}+\frac{t^2}{y}\right)
     =\partial_y^2+y^{-1}\partial_y+t^2y^{-2},
\]
for which
\begin{equation*}
 P_tK_{it}(2\pi n y)=(2\pi n)^2K_{it}(2\pi n y).
\end{equation*}
It is formally symmetric with respect to
\[
 \langle f,g\rangle=\int_0^\infty f(y)\overline{g(y)}\,y\,dy.
\]
With respect to this inner product, $P_t$ is formally self-adjoint: if
$f,g$ are smooth and at least one of them is compactly supported in
$(0,\infty)$, integration by parts gives
\[
 \langle P_tf,g\rangle
 =-\int_0^\infty y f'(y)\overline{g'(y)}\,dy
  +t^2\int_0^\infty f(y)\overline{g(y)}\,\frac{dy}{y}
 =\langle f,P_tg\rangle.
\]

\subsection{A uniform bound for the \texorpdfstring{$K$}{K}-Bessel function}

We use the following uniform estimate throughout this section.

\begin{lem}\label{lem:K-uniform}
Uniformly for $t\ge1$ and $x>0$,
\begin{equation*}
 e^{\pi t/2}|K_{it}(x)|\ll t^{-1/3}.
\end{equation*}
\end{lem}

\begin{proof}
For $x>1$, this is \cite[Proposition~4.3]{Liu2017}.  For $0<x\le1$, we use
\begin{equation}\label{eq:K-via-I}
 K_{it}(x)=\frac{\pi}{2i\sinh(\pi t)}
 \bigl(I_{-it}(x)-I_{it}(x)\bigr),
\end{equation}
which is formula~\cite[8.485]{GradshteynRyzhik2007}. We also have the power series for $I_{it}(x)$ \cite[8.445]{GradshteynRyzhik2007},
\[
 I_{it}(x)=\left(\frac{x}{2}\right)^{it}
 \sum_{k=0}^{\infty}
 \frac{(x^2/4)^k}{k!\,\Gamma(k+1+it)}.
\]
Since
\[
 |\Gamma(k+1+it)|
 =|\Gamma(1+it)|\prod_{j=1}^k|j+it|
 \ge |\Gamma(1+it)|\,k!,
\]
we obtain, uniformly for $0<x\le1$,
\[
 |I_{it}(x)|+|I_{-it}(x)|\ll |\Gamma(1+it)|^{-1}.
\]
Using $|\Gamma(1+it)|^2=\pi t/\sinh(\pi t)$ in
\eqref{eq:K-via-I} gives
\[
 e^{\pi t/2}|K_{it}(x)|\ll t^{-1/2}\ll t^{-1/3}.
\]
\end{proof}
We also need bounds for derivatives with respect to $x$. Since the normalising factor $e^{\pi t/2}$ is independent of $x$, the differential equation \eqref{eq:K-bessel-ode} also holds for $e^{\pi t/2}K_{it}(x).$

\begin{lem}\label{lem:K-euler-derivatives}
Fix $x_0>0$ and $M>0$.  For every integer $\ell\ge0$, uniformly for
$t\ge1$ and $x_0\le x\le Mt$,
\begin{equation}\label{eq:K-euler-derivatives}
 L_x^\ell\bigl(e^{\pi t/2}K_{it}(x)\bigr)\ll_{\ell,x_0,M}t^{\ell-1/3}.
\end{equation}
\end{lem}

\begin{proof}
We first treat $\ell=1$.  Write $x=e^s$ and
$u(s)=e^{\pi t/2}K_{it}(e^s)$.  Then
\[
 u'(s)=L_x\bigl(e^{\pi t/2}K_{it}(x)\bigr),\qquad
 u''(s)=(e^{2s}-t^2)u(s).
\]
With $h=t^{-1}$, Taylor's formula with integral remainder gives
\[
 u'(s)=\frac{u(s+h)-u(s)}{h}
 -\frac1h\int_0^h(h-r)u''(s+r)\,dr.
\]
For $0\le r\le h$ we have
$x_0\le e^{s+r}\le eMt$.  Hence Lemma~\ref{lem:K-uniform} yields
$u(s+r)\ll t^{-1/3}$, while
$u''(s+r)\ll t^{5/3}$.  Consequently
\[
 L_x\bigl(e^{\pi t/2}K_{it}(x)\bigr)=u'(s)\ll t^{2/3}.
\]

For the higher derivatives, we claim that for every $\ell\ge0$ there are
polynomials $\omega_{\ell,0}(x,t)$ and $\omega_{\ell,1}(x,t)$ such that
\begin{equation}\label{eq:euler-reduction}
 L_x^\ell\bigl(e^{\pi t/2}K_{it}(x)\bigr)
 =\omega_{\ell,0}(x,t)e^{\pi t/2}K_{it}(x)
  +\omega_{\ell,1}(x,t)L_x\bigl(e^{\pi t/2}K_{it}(x)\bigr),
\end{equation}
where, for $x\le Mt$,
\[
 \omega_{\ell,0}(x,t)\ll_{\ell,M}t^\ell,
 \qquad
 \omega_{\ell,1}(x,t)\ll_{\ell,M}t^{\ell-1}\quad(\ell\ge1).
\]
Indeed, take $\omega_{0,0}=1$, $\omega_{0,1}=0$ and use the Bessel equation to obtain
recursively
\[
 \omega_{\ell+1,0}=L_x\omega_{\ell,0}+(x^2-t^2)\omega_{\ell,1},
 \qquad
 \omega_{\ell+1,1}=\omega_{\ell,0}+L_x\omega_{\ell,1}.
\]
The stated polynomial bounds follow by induction. The estimate
\eqref{eq:K-euler-derivatives} then follows from \eqref{eq:euler-reduction},
Lemma~\ref{lem:K-uniform}, and the bound for $L_x(e^{\pi t/2}K_{it})$.
\end{proof}

We also need ordinary derivatives. For every integer $\ell\ge1$,
\[
 x^\ell\frac{d^\ell}{dx^\ell}
 =L_x(L_x-1)\cdots(L_x-\ell+1).
\]
Indeed, this follows by induction from
\[
 x^{\ell+1}\frac{d^{\ell+1}}{dx^{\ell+1}}
 =(L_x-\ell)x^\ell\frac{d^\ell}{dx^\ell}.
\]
Consequently, for suitable constants $c_{\ell,k}$,
\[
 \frac{d^\ell}{dx^\ell}
 =x^{-\ell}\sum_{k=0}^\ell c_{\ell,k}L_x^k.
\]
Since $x\ge x_0$, Lemma~\ref{lem:K-euler-derivatives} gives
\[
 e^{\pi t/2}\frac{d^\ell}{dx^\ell}K_{it}(x)
 \ll_{\ell,x_0,M}t^{\ell-1/3},
 \qquad x_0\le x\le Mt.
\]

For a positive integer $n$, set
\[
 F_n(y)=e^{\pi t/2}K_{it}(2\pi ny).
\]
If $L_y=y\partial_y$, then the scaling $x=2\pi ny$ gives
\begin{equation}\label{eq:euler-scaling}
 L_y^\ell F_n(y)=\left. L_x^\ell\bigl(e^{\pi t/2}K_{it}(x)\bigr)\right|_{x=2\pi ny}.
\end{equation}
Moreover,
\[
 y^j\partial_y^j=L_y(L_y-1)\cdots(L_y-j+1).
\]
Since $\supp\psi$ is contained in a fixed compact subinterval of
$(0,\infty)$, Lemma~\ref{lem:K-euler-derivatives} and
\eqref{eq:euler-scaling} imply that, for every fixed $j\ge0$,
\begin{equation}\label{eq:K-y-derivatives}
 \partial_y^jF_n(y)\ll_{\psi,j}t^{j-1/3},
 \qquad n\ll_\psi t,\quad y\in\supp\psi.
\end{equation}

\subsection{The Sturm--Liouville truncation}

\begin{lem}\label{lem:TR-structure}
Let
\[
 T_0=W,\qquad T_{R+1}=[P_t,T_R]=P_tT_R-T_RP_t,
 \qquad W(y)=\frac{\psi(y)}{y}.
\]
For every integer $R\ge0$ one has
\begin{equation}\label{eq:TR-structure}
 T_R=\sum_{j=0}^R c_{R,j}(y;t)\partial_y^j,
\end{equation}
where every coefficient $c_{R,j}$ is supported in $\supp\psi$.  Moreover,
for every integer $\ell\ge0$,
\begin{equation}\label{eq:TR-coefficients}
 \partial_y^\ell c_{R,j}(y;t)
 \ll_{\psi,R,\ell}t^{R-j}.
\end{equation}
\end{lem}

\begin{proof}
We use induction on $R$ to prove simultaneously the representation
\eqref{eq:TR-structure}, the support condition
\[
 \supp c_{R,j}\subseteq\supp\psi,
\]
and the estimates \eqref{eq:TR-coefficients} for every $\ell\ge0$.
For $R=0$ these assertions are immediate, since $T_0=W=\psi/y$ is
multiplication by a fixed smooth function supported in $\supp\psi$.

Suppose that they hold for $T_R$.  Then
\[
 T_{R+1}
 =\sum_{j=0}^R[P_t,c_{R,j}(y;t)\partial_y^j].
\]
Every coefficient occurring in these commutators is a linear combination
of derivatives of $c_{R,j}$ multiplied by smooth functions of $y$, and is
therefore supported in $\supp\psi$.  We estimate the contributions
from the three terms in $P_t$ separately.

First, Leibniz' rule gives
\begin{equation*}
 [\partial_y^2,c_{R,j}\partial_y^j]
 =2c_{R,j}'\partial_y^{j+1}+c_{R,j}''\partial_y^j.
\end{equation*}
Thus the first term contributes to $c_{R+1,j+1}$ and the second to
$c_{R+1,j}$.  By the induction hypothesis, for every $\ell\ge0$,
\[
 \partial_y^\ell c_{R,j}'\ll_{\psi,R,\ell}t^{R-j}
 =t^{R+1-(j+1)},
\]
while
\[
 \partial_y^\ell c_{R,j}''\ll_{\psi,R,\ell}t^{R-j}
 \le t^{R+1-j}.
\]
Both contributions therefore satisfy the required coefficient bounds at
their respective differential orders.

For the first-order term in $P_t$, a direct calculation gives
\begin{equation*}
 \left[\frac1y\partial_y,c_{R,j}\partial_y^j\right]
 =\frac{c_{R,j}'}y\partial_y^j
 -c_{R,j}\sum_{\nu=1}^j(-1)^\nu\nu!\binom{j}{\nu}
 y^{-\nu-1}\partial_y^{j-\nu+1}.
\end{equation*}
The first term contributes to $c_{R+1,j}$.  A term indexed by $\nu$ in
the sum contributes to $c_{R+1,k}$ with
\[
 k=j-\nu+1\le j.
\]
By Leibniz' rule and the induction hypothesis, for every $\ell\ge0$ the
$\ell$-th derivative of each of these coefficient contributions is
\[
 \ll_{\psi,R,\ell}t^{R-j}
 \le t^{R+1-k}.
\]

Finally,
\begin{equation*}
 \left[\frac{t^2}{y^2},c_{R,j}\partial_y^j\right]
 =-c_{R,j}\sum_{\nu=1}^j(-1)^\nu(\nu+1)!\binom{j}{\nu}
 t^2y^{-\nu-2}\partial_y^{j-\nu}.
\end{equation*}
A term indexed by $\nu$ contributes to $c_{R+1,k}$ with
\[
 k=j-\nu\le j-1.
\]
Again by Leibniz' rule and the induction hypothesis, for every
$\ell\ge0$ the $\ell$-th derivative of this contribution is
\[
 \ll_{\psi,R,\ell}t^{R-j+2}.
\]
Since $k\le j-1$,
\[
 R-j+2\le R+1-k,
\]
so this is $O_{\psi,R,\ell}(t^{R+1-k})$.

The three commutators also show that $T_{R+1}$ has differential order at
most $R+1$.
For each fixed $k$, summing the finitely many contributions to the
coefficient of $\partial_y^k$ gives
\[
 \partial_y^\ell c_{R+1,k}(y;t)
 \ll_{\psi,R,\ell}t^{R+1-k},
\]
for every $\ell\ge0$.  Together with the support property established
above, this completes the induction.
\end{proof}

\begin{lem}\label{lem:commutator}
Let $\psi\in C_c^\infty(0,\infty)$ and $m,n\ll_\psi t$.  For every integer
$R\ge0$,
\begin{equation}\label{eq:commutator}
 e^{\pi t}\left|\int_0^\infty
 \psi(y)K_{it}(2\pi my)K_{it}(2\pi ny)\,dy\right|
 \ll_{\psi,R}t^{-2/3}
 \left(1+\frac{|m^2-n^2|}{t}\right)^{-R}.
\end{equation}
\end{lem}

\begin{proof}
Put $\lambda_n=(2\pi n)^2$ and let $T_R$ be as in
Lemma~\ref{lem:TR-structure}.  Since
\[
 P_tF_n=\lambda_nF_n
\]
and the coefficients of $T_r$ are supported in $\supp\psi$, the formal
symmetry of $P_t$ gives
\[
\begin{aligned}
 \langle T_{r+1}F_n,F_m\rangle
 &=\langle P_tT_rF_n,F_m\rangle
   -\langle T_rP_tF_n,F_m\rangle\\
 &=\langle T_rF_n,P_tF_m\rangle
   -\lambda_n\langle T_rF_n,F_m\rangle\\
 &=(\lambda_m-\lambda_n)\langle T_rF_n,F_m\rangle.
\end{aligned}
\]
Iterating this identity yields
\begin{equation}\label{eq:commutator-exact}
 (\lambda_m-\lambda_n)^R\langle WF_n,F_m\rangle
 =\langle T_RF_n,F_m\rangle,
 \qquad W(y)=\frac{\psi(y)}{y}.
\end{equation}

By Lemma~\ref{lem:TR-structure} and \eqref{eq:K-y-derivatives},
\[
 |T_RF_n(y)|
 \ll_{\psi,R}\sum_{j=0}^Rt^{R-j}t^{j-1/3}
 \ll_{\psi,R}t^{R-1/3}
 \qquad (y\in\supp\psi).
\]
Together with Lemma~\ref{lem:K-uniform}, this gives
\begin{equation}\label{eq:TR-pairing}
 |\langle T_RF_n,F_m\rangle|
 \ll_{\psi,R}t^{R-2/3}.
\end{equation}
Hence, for $m\ne n$, \eqref{eq:commutator-exact} and
\eqref{eq:TR-pairing} imply
\[
 |\langle WF_n,F_m\rangle|
 \ll_{\psi,R}t^{-2/3}
 \left(\frac{t}{|m^2-n^2|}\right)^R.
\]
On the other hand, Lemma~\ref{lem:K-uniform} gives
\[
 |\langle WF_n,F_m\rangle|\ll_\psi t^{-2/3}.
\]
Combining these two estimates, and using
$|\lambda_m-\lambda_n|=(2\pi)^2|m^2-n^2|$, we obtain
\[
 |\langle WF_n,F_m\rangle|
 \ll_{\psi,R}t^{-2/3}
 \left(1+\frac{|m^2-n^2|}{t}\right)^{-R}.
\]
This proves \eqref{eq:commutator}.
\end{proof}

\subsection{The truncation error on a spectral interval}

Define
\[
 \mathcal E_j^{\mathrm{shift}}
 =\frac{1}{2L_j}
   \sum_{\substack{m\ne n\\m,n\le C_\psi T\\
                    |m^2-n^2|>T^{1+\eps}}}
   \lambda_j(m)\lambda_j(n)\cB_{m,n}(t_j).
\]

\begin{lem}
\label{lem:weighted-truncation}
Fix $\eps>0$ in the definition of $\mathcal E_j^{\mathrm{shift}}$.
For every $A>0$, uniformly for $1\le U\le T$,
\begin{equation*}
 \sum_{\substack{j\ {\rm even}\\T<t_j\le T+U}}L_j
 |\mathcal E_j^{\mathrm{shift}}|^2
 \ll_{A,\psi,\eps}T^{-A}.
\end{equation*}
\end{lem}

\begin{proof}
This is an immediate application of Lemma~\ref{lem:commutator}.
\end{proof}

By Lemmas~\ref{lem:K-tail} and \ref{lem:weighted-truncation}, we may
restrict the off-diagonal sum to
\[
 m,n\le C_\psi T,\qquad |m^2-n^2|\le T^{1+\eps}.
\]
Applying \eqref{eq:hecke} and writing
$m=da$, $n=d(a+b)$ gives
\begin{equation}\label{eq:trunc-main}
 d^2b(2a+b)=|d^2(a+b)^2-d^2a^2|\ll T^{1+\eps}.
\end{equation}

\section{Bounds for the Bessel integral and its \texorpdfstring{\(t\)}{t}-derivatives}\label{sec:bessel-bounds}

We use the truncation established in
Section~\ref{sec:trunc}, with the same parameter $\eps$. Put
\[
 p=da,\qquad q=d(a+b),
\]
so that
\[
 q^2-p^2=d^2b(2a+b)\ll T^{1+\eps}.
\]

Our goal is the uniform spectral derivative estimate in
Lemma~\ref{lem:B-derivatives}, with $\Lambda$ defined there:
\[
 \bigl|\partial_t^r\cB_{p,q}(t)\bigr|
 \ll_{r,\psi,\eps}\frac{\log T}{T}
 \Lambda\!\left(\frac{q-p}{p},T\right)^r.
\]
We begin by establishing the Fourier representation
\begin{equation}\label{eq:B-Fourier-preview}
 \cB_{p,q}(t)
 =2(1+e^{-2\pi t})
 \int_{\mathbb R}e^{it\xi}\mathscr A_{p,q}(\xi)\,d\xi,
\end{equation}
where, for $\xi\ne\pm\log(q/p)$, we define
\begin{equation*}
 \mathscr A_{p,q}(\xi)
 :=\int_{\mathbb R}
 \widehat\psi\bigl(p\sinh u+q\sinh(\xi-u)\bigr)\,du
\end{equation*}
 with the convention
\[\widehat\psi(\zeta)=\int_{\mathbb R}\psi(y)e^{-2\pi iy\zeta}\,dy.
\]

Then the analysis rests on the kernel $\mathscr A_{p,q}(\xi)$.
Its integral
representations change at the transition points
\[
 \xi=\pm\log(q/p),
\]
which separate the regimes of the kernel. We study its behaviour near
these points as well as in their complement, before turning to the
spectral derivatives of $\cB_{p,q}(t)$.

\subsection{The function \texorpdfstring{\(\mathscr A_{p,q}\)}{A(p,q)}}

We begin with a change of variables that describes the structure of
\(\mathscr A_{p,q}\).

\begin{lem}\label{lem:A-exact}
Let \(q>p\), and put
\[
 \vartheta=\log(q/p),\qquad
 \mathcal D_{p,q}(\xi):=p^2+q^2-2pq\cosh\xi.
\]
Then
\begin{align}
 \mathscr A_{p,q}(\xi)
 &=\int_{\mathbb R}\frac{\widehat\psi(s)}{\sqrt{s^2+\mathcal D_{p,q}(\xi)}}\,ds,
 && |\xi|<\vartheta,                                      \label{eq:A-interior}\\
 \mathscr A_{p,q}(\xi)
 &=2\int_{\sqrt{-\mathcal D_{p,q}(\xi)}}^\infty
   \frac{\widehat\psi(s)}{\sqrt{s^2+\mathcal D_{p,q}(\xi)}}\,ds,
 && \xi>\vartheta,                                        \label{eq:A-right}\\
 \mathscr A_{p,q}(\xi)
 &=2\int_{-\infty}^{-\sqrt{-\mathcal D_{p,q}(\xi)}}
   \frac{\widehat\psi(s)}{\sqrt{s^2+\mathcal D_{p,q}(\xi)}}\,ds,
 && \xi<-\vartheta.                                       \label{eq:A-left}
\end{align}
The square roots in \eqref{eq:A-right}--\eqref{eq:A-left} are positive on
the domains of integration.
\end{lem}

\begin{proof}
Put
\[
 F_\xi(u)=p\sinh u+q\sinh(\xi-u).
\]
We can write
\[
 F_\xi(u)=\frac12\{\alpha(\xi)e^u+\beta(\xi)e^{-u}\},
 \qquad
 \alpha(\xi)=p-qe^{-\xi},\quad
 \beta(\xi)=qe^\xi-p.
\]
Then
\[
 \partial_uF_\xi(u)
 =\frac12\{\alpha(\xi)e^u-\beta(\xi)e^{-u}\},
\]
and hence
\[
 \{\partial_uF_\xi(u)\}^2-F_\xi(u)^2
 =-\alpha(\xi)\beta(\xi)
 =p^2+q^2-2pq\cosh\xi
 =\mathcal D_{p,q}(\xi).
\]
Thus
\begin{equation}\label{eq:F-invariant}
 \{\partial_uF_\xi(u)\}^2-F_\xi(u)^2=\mathcal D_{p,q}(\xi).
\end{equation}
Also,
\begin{equation*}
 \mathcal D_{p,q}(\xi)
 =4pq\sinh\frac{\vartheta+\xi}{2}\sinh\frac{\vartheta-\xi}{2}.
\end{equation*}
Thus \(\mathcal D_{p,q}\) has simple zeros at \(\xi=\pm\vartheta\), with
\[
 \mathcal D_{p,q}'(\vartheta)=-(q^2-p^2),\qquad
 \mathcal D_{p,q}'(-\vartheta)=q^2-p^2.
\]
The points \(\xi=\pm\vartheta\) separate the region where \(F_\xi\) is monotone from the regions where it has a critical point.

Suppose first that $|\xi|<\vartheta$.  Then
\[
 \alpha(\xi)<0,\qquad \beta(\xi)>0,\qquad
 \alpha(\xi)\beta(\xi)<0.
\]
By \eqref{eq:F-invariant},
\[
 |\partial_uF_\xi(u)|
 =\sqrt{F_\xi(u)^2+\mathcal D_{p,q}(\xi)}>0.
\]
Moreover,
\[
 \lim_{u\to-\infty}F_\xi(u)=+\infty,\qquad
 \lim_{u\to+\infty}F_\xi(u)=-\infty.
\]
Thus $F_\xi$ is strictly decreasing and maps $\mathbb R$ bijectively onto
$\mathbb R$.  With $s=F_\xi(u)$,
\begin{equation*}
 |du|
 =\frac{ds}{|\partial_uF_\xi(u)|}
 =\frac{ds}{\sqrt{s^2+\mathcal D_{p,q}(\xi)}}.
\end{equation*}
Substituting this into the integral defining \(\mathscr A_{p,q}(\xi)\) gives \eqref{eq:A-interior}.

If $\xi>\vartheta$, then $\alpha(\xi),\beta(\xi)>0$.  The function
$F_\xi$ is positive and strictly convex, with a unique minimum at
\[
 e^{2u_0}=\frac{\beta(\xi)}{\alpha(\xi)}.
\]
At this point,
\[
 F_\xi(u_0)
 =\sqrt{\alpha(\xi)\beta(\xi)}
 =\sqrt{-\mathcal D_{p,q}(\xi)}.
\]
Hence the range is
$[\sqrt{-\mathcal D_{p,q}(\xi)},\infty)$, and every value above the minimum
has two preimages.  On both branches the absolute Jacobian equals
$(s^2+\mathcal D_{p,q}(\xi))^{-1/2}$, which gives
\eqref{eq:A-right}.

If $\xi<-\vartheta$, then $\alpha(\xi),\beta(\xi)<0$.  The function is
strictly concave, has maximum
$-\sqrt{-\mathcal D_{p,q}(\xi)}$, and has range
$(-\infty,-\sqrt{-\mathcal D_{p,q}(\xi)}]$.  The same two-branch change of
variables proves \eqref{eq:A-left}.  All square roots are taken positive.
\end{proof}

The representations in Lemma~\ref{lem:A-exact} reduce the behaviour of $\mathscr A_{p,q}$ near the zeros of $\mathcal D_{p,q}$ to the following model integrals.

\begin{lem}
\label{lem:model-integrals}
Let $\varphi\in\mathcal S(\mathbb R)$.  For every $h>0$, define
\begin{align*}
 \Phi_0(h)&=\int_{\mathbb R}\frac{\varphi(s)}{\sqrt{s^2+h}}\,ds,\\
 \Phi_+(h)&=2\int_{\sqrt h}^{\infty}
       \frac{\varphi(s)}{\sqrt{s^2-h}}\,ds,\\
 \Phi_-(h)&=2\int_{-\infty}^{-\sqrt h}
       \frac{\varphi(s)}{\sqrt{s^2-h}}\,ds.
\end{align*}
For every $\nu\in\{0,+,-\}$ there are a constant $c_\nu(\varphi)$ and a
function $R_\nu\in C^\infty((0,\infty))$ such that
\begin{equation}\label{eq:model-log}
 \Phi_\nu(h)=-\varphi(0)\log h+c_\nu(\varphi)+R_\nu(h)
\end{equation}
and, for $0<h\le1$ and every fixed $\ell\ge0$,
\begin{equation}\label{eq:model-remainder-derivatives}
 |R_\nu^{(\ell)}(h)|\ll_{\ell,\varphi}
 h^{1/2-\ell}(1+|\log h|).
\end{equation}
For $h\ge1$, $\Phi_\pm$ and all their fixed $h$-derivatives are
$O_{N,\ell,\varphi}(h^{-N})$.  If all moments
$\int_{\mathbb R}s^k\varphi(s)\,ds$ vanish, the same assertion holds for
$\Phi_0$.
\end{lem}

\begin{proof}
We first prove \eqref{eq:model-log}--\eqref{eq:model-remainder-derivatives}
for $0<h\le1/2$.

We begin with
\[
 \Phi_0(h)=\int_{\mathbb R}\frac{\varphi(s)}{\sqrt{s^2+h}}\,ds.
\]
Choose an even $\chi\in C_c^\infty((-2,2))$ with $\chi=1$ on $[-1,1]$, and put
\[
 g(s)=\varphi(s)-\varphi(0)\chi(s).
\]
Then $g(s)=O_\varphi(|s|)$ near zero and
\[
 \Phi_0(h)=\varphi(0)\int_{\mathbb R}\frac{\chi(s)}{\sqrt{s^2+h}}\,ds
 +\int_{\mathbb R}\frac{g(s)}{\sqrt{s^2+h}}\,ds.
\]
Since $\chi=1$ on $[-1,1]$,
\begin{align*}
 \int_{\mathbb R}\frac{\chi(s)}{\sqrt{s^2+h}}\,ds
 &=2\int_0^1\frac{ds}{\sqrt{s^2+h}}
   +\int_{1\le |s|\le2}\frac{\chi(s)}{\sqrt{s^2+h}}\,ds\\
 &=-\log h+2\log2
   +\int_{1\le |s|\le2}\frac{\chi(s)}{|s|}\,ds\\
 &\quad+2\log\frac{1+\sqrt{1+h}}2
   +\int_{1\le |s|\le2}\chi(s)
   \left(\frac1{\sqrt{s^2+h}}-\frac1{|s|}\right)ds,
\end{align*}
where
\[
 2\int_0^1\frac{ds}{\sqrt{s^2+h}}
 =2\operatorname{arsinh}(h^{-1/2})
 =-\log h+2\log(1+\sqrt{1+h}).
\]
Define
\[
 c_0(\varphi)
 :=\varphi(0)\left(2\log2+
 \int_{1\le |s|\le2}\frac{\chi(s)}{|s|}\,ds\right)
 +\int_{\mathbb R}\frac{g(s)}{|s|}\,ds
\]
and
\begin{align*}
 R_0(h)
 &:=2\varphi(0)\log\frac{1+\sqrt{1+h}}2\\
 &\quad+\varphi(0)\int_{1\le |s|\le2}\chi(s)
 \left(\frac1{\sqrt{s^2+h}}-\frac1{|s|}\right)ds\\
 &\quad+\int_{\mathbb R}g(s)
 \left(\frac1{\sqrt{s^2+h}}-\frac1{|s|}\right)ds.
\end{align*}
Thus
\[
 \Phi_0(h)=-\varphi(0)\log h+c_0(\varphi)+R_0(h).
\]
For $0<h\le1/2$,
\[
 \log\frac{1+\sqrt{1+h}}2\ll h,
 \qquad
 \partial_h^\ell\log\frac{1+\sqrt{1+h}}2\ll_\ell1
 \quad(\ell\ge1),
\]
and
\begin{align*}
 &\int_{1\le |s|\le2}\chi(s)
 \left(\frac1{\sqrt{s^2+h}}-\frac1{|s|}\right)ds\ll_\chi h,\\
 &\partial_h^\ell\int_{1\le |s|\le2}\chi(s)
 \left(\frac1{\sqrt{s^2+h}}-\frac1{|s|}\right)ds
 \ll_{\ell,\chi}1\qquad(\ell\ge1).
\end{align*}
For the last term, splitting at $|s|=\sqrt h$ and using
\[
 \left|\frac1{\sqrt{s^2+h}}-\frac1{|s|}\right|
 \ll
 \begin{cases}
 |s|^{-1},&|s|\le\sqrt h,\\
 h|s|^{-3},&|s|>\sqrt h,
 \end{cases}
\]
gives
\[
 \int_{\mathbb R}|g(s)|
 \left|\frac1{\sqrt{s^2+h}}-\frac1{|s|}\right|ds
 \ll_\varphi h^{1/2}.
\]
For every $\ell\ge1$, differentiation under the integral sign gives
\[
 \partial_h^\ell\int_{\mathbb R}\frac{g(s)}{\sqrt{s^2+h}}\,ds
 =c_\ell\int_{\mathbb R}g(s)(s^2+h)^{-\ell-1/2}\,ds
 \ll_{\ell,\varphi}h^{1/2-\ell},
\]
where $c_\ell$ depends only on $\ell$.  Hence
\[
 |R_0^{(\ell)}(h)|\ll_{\ell,\varphi}
 h^{1/2-\ell}(1+|\log h|),
 \qquad 0<h\le1/2.
\]

We next consider
\[
 \Phi_+(h)=2\int_{\sqrt h}^{\infty}
 \frac{\varphi(s)}{\sqrt{s^2-h}}\,ds.
\]
On $0\le s\le1$, put
\[
 \varphi_+(s)=\begin{cases}
 \dfrac{\varphi(s)-\varphi(0)}s,&s>0,\\[2mm]
 \varphi'(0),&s=0,
 \end{cases}
\]
so that $\varphi_+\in C^\infty([0,1])$ and
$\varphi(s)=\varphi(0)+s\varphi_+(s)$.  Splitting at $s=1$ and using
$u=\sqrt{s^2-h}$ gives
\begin{align*}
 \Phi_+(h)
 &=2\varphi(0)\operatorname{arccosh}(h^{-1/2})
   +2\int_0^{\sqrt{1-h}}
      \varphi_+\!\left(\sqrt{u^2+h}\right)du\\
 &\quad+2\int_1^\infty\frac{\varphi(s)}{\sqrt{s^2-h}}\,ds.
\end{align*}
Since
\[
 2\operatorname{arccosh}(h^{-1/2})
 =-\log h+2\log2+2\log\frac{1+\sqrt{1-h}}2,
\]
define
\[
 c_+(\varphi)
 :=2\varphi(0)\log2
   +2\int_0^1\varphi_+(u)du
   +2\int_1^\infty\frac{\varphi(s)}s\,ds
\]
and
\begin{align*}
 R_+(h)
 &:=2\varphi(0)\log\frac{1+\sqrt{1-h}}2\\
 &\quad+2\left\{
   \int_0^{\sqrt{1-h}}\varphi_+\!\left(\sqrt{u^2+h}\right)du
   -\int_0^1\varphi_+(u)du\right\}\\
 &\quad+2\int_1^\infty\varphi(s)
   \left(\frac1{\sqrt{s^2-h}}-\frac1s\right)ds.
\end{align*}
Then
\[
 \Phi_+(h)=-\varphi(0)\log h+c_+(\varphi)+R_+(h).
\]
For $0<h\le1/2$,
\[
 \log\frac{1+\sqrt{1-h}}2\ll h,
 \qquad
 \partial_h^\ell\log\frac{1+\sqrt{1-h}}2\ll_\ell1
 \quad(\ell\ge1),
\]
and
\begin{align*}
 &\int_1^\infty|\varphi(s)|
 \left|\frac1{\sqrt{s^2-h}}-\frac1s\right|ds
 \ll_\varphi h,\\
 &\partial_h^\ell\int_1^\infty\varphi(s)
 \left(\frac1{\sqrt{s^2-h}}-\frac1s\right)ds
 \ll_{\ell,\varphi}1\qquad(\ell\ge1).
\end{align*}
Moreover,
\[
 \int_0^{\sqrt{1-h}}
 \left|\varphi_+\!\left(\sqrt{u^2+h}\right)-\varphi_+(u)\right|du
 +\int_{\sqrt{1-h}}^1|\varphi_+(u)|du
 \ll_\varphi h(1+|\log h|).
\]
By Theorem~\ref{thm:faa-di-bruno}, for $\ell\ge1$,
\begin{equation*}
 \partial_h^\ell\varphi_+\!\left(\sqrt{u^2+h}\right)
 =\sum_{k=1}^\ell c_{\ell,k}\,
 \varphi_+^{(k)}\!\left(\sqrt{u^2+h}\right)
 (u^2+h)^{k/2-\ell},
\end{equation*}
where the constants $c_{\ell,k}$ depend only on $\ell$ and $k$.
The terms from the moving endpoint are bounded, since $u^2+h=1$ there, and
\[
 \partial_h^\ell
 \int_0^{\sqrt{1-h}}
 \varphi_+\!\left(\sqrt{u^2+h}\right)du
 \ll_{\ell,\varphi}
 \begin{cases}
 1+|\log h|,&\ell=1,\\
 h^{1-\ell},&\ell\ge2.
 \end{cases}
\]
Consequently,
\[
 |R_+^{(\ell)}(h)|\ll_{\ell,\varphi}
 h^{1/2-\ell}(1+|\log h|),
 \qquad 0<h\le1/2.
\]
Replacing $\varphi(s)$ by $\varphi(-s)$ gives the same conclusion for
$R_-(h)$.  Since $R_0,R_+,R_-\in C^\infty((0,\infty))$, their fixed
derivatives are bounded on $[1/2,1]$.  Thus
\eqref{eq:model-remainder-derivatives} holds for $0<h\le1$.

It remains to consider $h\ge1$.  We have
\[
 \Phi_\pm(h)=2\int_0^\infty
 \varphi(\pm\sqrt h\cosh v)\,dv.
\]
For $\ell\ge1$, Theorem~\ref{thm:faa-di-bruno} gives
\begin{equation*}
 \partial_h^\ell\varphi(\pm\sqrt h\cosh v)
 =\sum_{k=1}^\ell c_{\ell,k}^{\pm}\,
 h^{k/2-\ell}(\cosh v)^k
 \varphi^{(k)}(\pm\sqrt h\cosh v).
\end{equation*}
Together with the case $\ell=0$, the Schwartz decay of $\varphi$ yields
\[
 \Phi_\pm^{(\ell)}(h)\ll_{N,\ell,\varphi}h^{-N}
\]
for every $N>0$.

Finally,
\[
 \Phi_0^{(\ell)}(h)
 =c_\ell\int_{\mathbb R}
 \varphi(s)(s^2+h)^{-\ell-1/2}\,ds.
\]
Assume that all moments of $\varphi$ vanish.  On $|s|\le\sqrt{h/2}$,
expand the last factor to an arbitrarily high order in $s^2/h$.  For every
$k\ge0$,
\[
 \int_{|s|\le\sqrt{h/2}}s^{2k}\varphi(s)\,ds
 =-\int_{|s|>\sqrt{h/2}}s^{2k}\varphi(s)\,ds
 \ll_{N,k,\varphi}h^{-N}.
\]
The expansion remainder and the range $|s|>\sqrt{h/2}$ are bounded in the
same way by the Schwartz decay of $\varphi$.  Therefore
\[
 \Phi_0^{(\ell)}(h)\ll_{N,\ell,\varphi}h^{-N}
\]
for every $N>0$, completing the proof.
\end{proof}

In our application $\varphi=\widehat\psi$.  Since
$\supp\psi\Subset(0,\infty)$, for every $k\ge0$,
\begin{equation*}
 \int_{\mathbb R}s^k\widehat\psi(s)\,ds
 =(2\pi i)^{-k}\psi^{(k)}(0)=0.
\end{equation*}

We first recall the definition of the total variation used in the next lemma.
If
$I=[\alpha,\beta]$ and $F:I\to\mathbb C$, set
\begin{equation*}
 \operatorname{Var}_{I}(F)
 :=\sup_{\alpha=x_0<\cdots<x_m=\beta}
 \sum_{k=1}^{m}|F(x_k)-F(x_{k-1})|,
\end{equation*}
where the supremum is over all finite partitions of $I$; finite total
variation means $F\in BV(I)$. Moreover, if $F$ is absolutely continuous, then
\begin{equation*}
 \operatorname{Var}_{I}(F)=\int_I|F'(x)|\,dx.
\end{equation*}
If $h$ is monotone on $I$ and $G$ is absolutely continuous on
$h(I)$, then
\begin{equation}\label{eq:variation-composition}
 \operatorname{Var}_{I}(G\circ h)
 \le \int_{h(I)}|G'(u)|\,du.
\end{equation}

\begin{lem}\label{lem:A-composition}
Assume $q>p\ge1$, $p,q\ll T$, and
$1\le q^2-p^2\ll T^{1+\eps}$.  For every fixed $r\ge0$,
\begin{equation}\label{eq:A-absolute}
 \int_{\mathbb R}|\xi|^r|\mathscr A_{p,q}(\xi)|\,d\xi<\infty.
\end{equation}
Let $0<\Delta\le1/4$.  For fixed $j,N\ge0$,
\begin{equation}\label{eq:A-derivative-away}
 \int_{\operatorname{dist}(\xi,\{\pm\vartheta\})\ge\Delta}
 (1+|\xi|)^N|\partial_\xi^j\mathscr A_{p,q}(\xi)|\,d\xi
 \ll_{j,N,\psi}(\log T)^N\Delta^{-j}.
\end{equation}
For $\sigma\in\{+,-\}$ put
\[
 I_\sigma=[\sigma\vartheta-2\Delta,\sigma\vartheta+2\Delta].
\]
There are explicitly constructed functions $H_\sigma\in BV(I_\sigma)$ such
that, for $\xi\in I_\sigma\setminus\{\pm\vartheta\}$,
\begin{equation}\label{eq:A-local-log}
 \mathscr A_{p,q}(\xi)=-\widehat\psi(0)\log|\mathcal D_{p,q}(\xi)|+H_\sigma(\xi),
\end{equation}
and
\begin{equation}\label{eq:H-BV}
 \|H_\sigma\|_{L^\infty(I_\sigma)}+
 \operatorname{Var}_{I_\sigma}(H_\sigma)
 \ll_\psi\log T.
\end{equation}
\end{lem}

\begin{proof}
Put
\[
 h(\xi):=|\mathcal D_{p,q}(\xi)|,
 \qquad
 \nu(\xi)=
 \begin{cases}
  -, & \xi<-\vartheta,\\
  0, & |\xi|<\vartheta,\\
  +, & \xi>\vartheta.
 \end{cases}
\]
We take $\varphi=\widehat\psi$ in Lemma~\ref{lem:model-integrals}.  By
Lemma~\ref{lem:A-exact},
\begin{equation}\label{eq:A-model-selector}
 \mathscr A_{p,q}(\xi)=\Phi_{\nu(\xi)}(h(\xi)),
 \qquad \xi\ne\pm\vartheta.
\end{equation}
Thus the proof follows the change of the single parameter $h(\xi)$ as
$\xi$ moves along the real line.  The zeros $\pm\vartheta$ divide the line
as indicated below; on each interval obtained by also splitting at $0$,
$h(\xi)$ is monotone.
\begin{center}
\begin{tikzpicture}[x=1.05cm,y=1cm]
 \draw[->] (-5.2,0) -- (5.2,0) node[right] {$\xi$};
 \foreach \x in {-3.2,-1.2,1.2,3.2}
   \draw (\x,0.09)--(\x,-0.09);
 \draw (-2.2,0.11)--(-2.2,-0.11) node[below=4pt] {$-\vartheta$};
 \draw (2.2,0.11)--(2.2,-0.11) node[below=4pt] {$\vartheta$};
 \draw[decorate,decoration={brace,amplitude=4pt}]
   (-3.2,0.42)--(-1.2,0.42) node[midway,above=6pt] {$4\Delta$};
 \draw[decorate,decoration={brace,amplitude=4pt}]
   (1.2,0.42)--(3.2,0.42) node[midway,above=6pt] {$4\Delta$};
\end{tikzpicture}
\end{center}

We first prove \eqref{eq:A-local-log} and \eqref{eq:H-BV}.  On $I_+\cup I_-$ we have $|\xi|\le\vartheta+1/2$ and $p,q\ll T$; hence
\begin{equation*}
 0\le h(\xi)\le CT^2
\end{equation*}
for a constant $C$ depending only on the fixed support data.  Thus, for these two assertions, it is enough to control the model functions on $0<h\le CT^2$.

For $\nu\in\{0,+,-\}$ define
\begin{equation}\label{eq:model-regular-part}
 F_\nu(h):=\Phi_\nu(h)+\widehat\psi(0)\log h.
\end{equation}
Lemma~\ref{lem:model-integrals} gives
\[
 F_\nu(h)=c_\nu(\widehat\psi)+R_\nu(h),\qquad 0<h\le1,
\]
and therefore
\[
 \|F_\nu\|_{L^\infty(0,1]}+
 \int_0^1|F_\nu'(u)|\,du\ll_\psi1.
\]
For $1\le h\le CT^2$, the large-parameter estimates in
Lemma~\ref{lem:model-integrals} give
\[
 |F_\nu(h)|\ll_\psi1+\log h,
 \qquad |F_\nu'(h)|\ll_\psi h^{-1}.
\]
Hence
\begin{equation}\label{eq:model-BV-range}
 \|F_\nu\|_{L^\infty(0,CT^2]}+
 \int_0^{CT^2}|F_\nu'(u)|\,du
 \ll_\psi\log T.
\end{equation}

Extend each $F_\nu$ continuously to $h=0$ by setting
$F_\nu(0)=c_\nu(\widehat\psi)$. For $\sigma\in\{+,-\}$, define
\begin{equation*}
 H_\sigma(\xi):=F_{\nu(\xi)}(h(\xi)),
 \qquad \xi\in I_\sigma\setminus\{\pm\vartheta\},
\end{equation*}
and at an interior zero choose its right-hand trace. At either endpoint
of $I_\sigma$ choose the trace from inside $I_\sigma$; this convention
also applies when an endpoint is a zero. These finite traces exist by
\eqref{eq:model-log}--\eqref{eq:model-remainder-derivatives}.
Thus $H_\sigma$ is specified on the whole closed interval; a jump at a
zero is allowed. Away from the zeros, \eqref{eq:A-model-selector} and
\eqref{eq:model-regular-part} give
\[
 \mathscr A_{p,q}(\xi)
 =-\widehat\psi(0)\log h(\xi)+H_\sigma(\xi)
 =-\widehat\psi(0)\log|\mathcal D_{p,q}(\xi)|+H_\sigma(\xi),
\]
which proves \eqref{eq:A-local-log}. This identity is only asserted
away from the zeros, or equivalently almost everywhere on $I_\sigma$.

We next follow $h(\xi)$ on $I_\sigma$.  Partition $I_\sigma$ at the points among
$-\vartheta,0,\vartheta$ that it contains.  On each resulting interval $J$,
$\nu(\xi)$ is fixed and $h(\xi)$ is monotone.  Therefore
\eqref{eq:variation-composition} and \eqref{eq:model-BV-range} give
\[
 \operatorname{Var}_{J}(H_\sigma)
 \le \int_0^{CT^2}|F_\nu'(u)|\,du
 \ll_\psi\log T.
\]
At $\xi=\vartheta$ or $-\vartheta$, Lemma~\ref{lem:model-integrals} gives
$F_\nu(h)\to c_\nu(\widehat\psi)$ as $h\downarrow0$.  Thus crossing a zero
adds only one of the finite jumps
\[
 |c_+(\widehat\psi)-c_0(\widehat\psi)|,
 \qquad
 |c_-(\widehat\psi)-c_0(\widehat\psi)|.
\]
There are only finitely many pieces, so \eqref{eq:H-BV} follows.  When $I_+$ and $I_-$ overlap, we split them at the points among
$-\vartheta,0,\vartheta$. On each resulting interval, $\nu$ is fixed
and $h$ is monotone, so the same variation estimate applies.

It remains to prove \eqref{eq:A-derivative-away}.  Throughout this part suppose
\[
 \operatorname{dist}(\xi,\{\pm\vartheta\})\ge\Delta.
\]
The factorisation of $\mathcal D_{p,q}$ gives
\[
 \log h(\xi)
 =\log(4pq)+\log\left|\sinh\frac{\vartheta+\xi}{2}\right|
 +\log\left|\sinh\frac{\vartheta-\xi}{2}\right|,
\]
and hence, for every fixed $\ell\ge1$,
\begin{equation}\label{eq:logD-derivatives}
 |\partial_\xi^\ell\log h(\xi)|\ll_\ell\Delta^{-\ell}.
\end{equation}
On each interval determined by $-\vartheta,0,\vartheta$, the function $h$ is
smooth.  Applying Theorem~\ref{thm:faa-di-bruno} to
$h=\exp(\log h)$ gives
\begin{equation}\label{eq:D-derivatives-relative}
 |h^{(\ell)}(\xi)|\ll_\ell h(\xi)\Delta^{-\ell}.
\end{equation}

We first treat $0<h(\xi)\le1$.  By Lemma~\ref{lem:model-integrals},
\[
 \mathscr A_{p,q}(\xi)
 =-\widehat\psi(0)\log h(\xi)
  +c_{\nu(\xi)}(\widehat\psi)+R_{\nu(\xi)}(h(\xi)).
\]
For $j\ge1$, the logarithmic term is controlled by \eqref{eq:logD-derivatives}.  For the
remainder, Theorem~\ref{thm:faa-di-bruno},
\eqref{eq:D-derivatives-relative}, and
\eqref{eq:model-remainder-derivatives} show that every term in
$\partial_\xi^jR_\nu(h(\xi))$ is bounded by
\[
 h(\xi)^{1/2-k}(1+|\log h(\xi)|)
 \prod_{\ell=1}^j
 \bigl(h(\xi)\Delta^{-\ell}\bigr)^{m_\ell}
 =h(\xi)^{1/2}(1+|\log h(\xi)|)\Delta^{-j}
 \ll\Delta^{-j},
\]
where $k=\sum_{\ell=1}^jm_\ell$ and
$\sum_{\ell=1}^j\ell m_\ell=j$.  Hence, for every fixed $j\ge1$,
\begin{equation}\label{eq:A-small-h-derivative}
 |\partial_\xi^j\mathscr A_{p,q}(\xi)|\ll_{j,\psi}\Delta^{-j},
 \qquad 0<h(\xi)\le1.
\end{equation}
For $j=0$, the small-$h$ expansion instead gives
\begin{equation}\label{eq:A-small-h-zero}
 |\mathscr A_{p,q}(\xi)|\ll_\psi 1+|\log h(\xi)|,
 \qquad 0<h(\xi)\le1,
\end{equation}
since $h^{1/2}(1+|\log h|)$ is bounded on $(0,1]$.

We next treat $h(\xi)\ge1$.  Using \eqref{eq:A-model-selector} directly,
the large-$h$ estimates of Lemma~\ref{lem:model-integrals},
\eqref{eq:D-derivatives-relative}, and Theorem~\ref{thm:faa-di-bruno} give,
for every fixed $A>0$,
\begin{equation}\label{eq:A-large-h-derivative}
 |\partial_\xi^j\mathscr A_{p,q}(\xi)|
 \ll_{A,j,\psi}\Delta^{-j}h(\xi)^{-A}.
\end{equation}

It remains to integrate these bounds.  We use the elementary inequality
\begin{equation}\label{eq:D-level-lower}
 h(\xi)=2pq|\cosh\xi-\cosh\vartheta|
 \ge pq|\xi^2-\vartheta^2|
 \ge |\xi^2-\vartheta^2|.
\end{equation}
Indeed, $|\cosh x-\cosh y|\ge |x^2-y^2|/2$ for
$x,y\ge0$.  For $u>0$, it follows that
\begin{equation}\label{eq:D-sublevel-measure}
 |\{\xi:h(\xi)<u\}|
 \le |\{\xi:|\xi^2-\vartheta^2|<u\}|
 \le 4\sqrt u.
\end{equation}
To see the last inequality, if $\vartheta^2\le u$ the latter
set has length at most $2\sqrt{2u}$; otherwise its length is
$2(\sqrt{\vartheta^2+u}-\sqrt{\vartheta^2-u})
\le2\sqrt{2u}$.  Moreover, $h(\xi)<u$ implies
$|\xi|\le\vartheta+\sqrt u$.

On $0<h(\xi)\le1$, these observations give
$|\xi|\le\vartheta+1$ and a set of measure at most $4$.
For $j\ge1$, \eqref{eq:A-small-h-derivative} therefore yields
\[
 \int_{\substack{0<h(\xi)\le1\\
             \operatorname{dist}(\xi,\{\pm\vartheta\})\ge\Delta}}
 (1+|\xi|)^N|\partial_\xi^j\mathscr A_{p,q}(\xi)|\,d\xi
 \ll_{j,N,\psi}(1+\vartheta)^N\Delta^{-j}.
\]
For $j=0$, the logarithmic term in \eqref{eq:A-small-h-zero}
is integrable uniformly in $p,q$.  Indeed, \eqref{eq:D-sublevel-measure} gives
\begin{equation}\label{eq:D-log-L1}
 \int_{\{0<h(\xi)\le1\}}|\log h(\xi)|\,d\xi
 =\int_0^\infty
 |\{\xi:0<h(\xi)<e^{-v}\}|\,dv
 \le4\int_0^\infty e^{-v/2}\,dv=8.
\end{equation}
Together with \eqref{eq:A-small-h-zero}, this proves the same
weighted integral bound for $j=0$ (with $\Delta^0=1$).

For $h(\xi)\ge1$, divide the domain into the sets
$\Omega_k=\{\xi:2^k\le h(\xi)<2^{k+1}\}$, $k\ge0$.
By \eqref{eq:D-sublevel-measure},
$|\Omega_k|\le4\sqrt{2^{k+1}}\ll2^{k/2}$.
Moreover, \eqref{eq:D-level-lower} gives, for $\xi\in\Omega_k$,
\[
 |\xi|^2\le\vartheta^2+h(\xi)
 <\vartheta^2+2^{k+1},
 \qquad
 |\xi|\le\vartheta+2^{(k+1)/2}.
\]
Consequently,
\[
 (1+|\xi|)^N
 \le(1+\vartheta+2^{(k+1)/2})^N
 \ll_N(1+\vartheta)^N2^{kN/2}.
\]
On $\Omega_k$, \eqref{eq:A-large-h-derivative} also gives
$|\partial_\xi^j\mathscr A_{p,q}(\xi)|
 \ll_{B,j,\psi}\Delta^{-j}2^{-kB}$.
Thus, for any fixed $B>(N+1)/2$,
\[
 \begin{aligned}
 &\int_{\substack{h(\xi)\ge1\\
             \operatorname{dist}(\xi,\{\pm\vartheta\})\ge\Delta}}
 (1+|\xi|)^N|\partial_\xi^j\mathscr A_{p,q}(\xi)|\,d\xi\\
 &\quad\ll_{B,j,N,\psi}
 (1+\vartheta)^N\Delta^{-j}
 \sum_{k\ge0}2^{kN/2}\,2^{k/2}\,2^{-kB}\\
 &\quad\ll_{j,N,\psi}(1+\vartheta)^N\Delta^{-j},
 \end{aligned}
\]
since the resulting geometric series converges.
Combining the two ranges and using
$\vartheta=\log(q/p)\ll_\psi\log T$ proves
\eqref{eq:A-derivative-away} with the stated exponent $N$.

Finally, near $\pm\vartheta$ the representation
\eqref{eq:A-local-log} is locally integrable, while for large $h(\xi)$ the
functions $\Phi_\nu(h(\xi))$ decay rapidly.  This proves
\eqref{eq:A-absolute}.  If the two zero-neighbourhoods overlap, we partition at the points
among $-\vartheta,0,\vartheta$. On each resulting interval, $\nu$ is
fixed and $h$ is monotone, so the same bound applies.
\end{proof}

\subsection{A Fourier representation of the Bessel integral}
\label{subsec:exact-fourier}

By the preceding lemmas, for every fixed integer $m\ge0$,
\begin{equation}\label{eq:A-weighted-L1}
 \int_{\mathbb R}(1+|\xi|)^m
 |\mathscr A_{p,q}(\xi)|\,d\xi<\infty.
\end{equation}
\begin{lem}
\label{lem:exact-fourier-t}
For $q>p\ge1$, recall
\[
 \cB_{p,q}(t):=16\cosh(\pi t)
 \int_0^\infty\psi(y)K_{it}(2\pi py)K_{it}(2\pi qy)\,dy.
\]
Then
\begin{equation}\label{eq:B-Fourier-exact}
 \cB_{p,q}(t)
 =2(1+e^{-2\pi t})
 \int_{\mathbb R}e^{it\xi}\mathscr A_{p,q}(\xi)\,d\xi.
\end{equation}
For every fixed integer $r\ge0$,
\begin{align}
 \partial_t^r\cB_{p,q}(t)
 &=2(1+e^{-2\pi t})
 \int_{\mathbb R}(i\xi)^re^{it\xi}
 \mathscr A_{p,q}(\xi)\,d\xi \notag\\
 &\quad
 +2e^{-2\pi t}\sum_{\ell=1}^r
 \binom r\ell(-2\pi)^\ell
 \int_{\mathbb R}(i\xi)^{r-\ell}e^{it\xi}
 \mathscr A_{p,q}(\xi)\,d\xi,                    \label{eq:B-derivative-exact}
\end{align}
with an empty second line when $r=0$.
\end{lem}

\begin{proof}
For $0<\eta<\pi/2$, shifting the standard integral representation
\[
 K_{it}(x)=\frac12\int_{\mathbb R}e^{-x\cosh z}e^{itz}\,dz
\]
to $\Im z=\pi/2-\eta$ gives
\begin{equation}\label{eq:shifted-K}
 2e^{(\pi/2-\eta)t}K_{it}(x)
 =\int_{\mathbb R}e^{-x\sin\eta\cosh u}
 e^{i\{tu-x\cos\eta\sinh u\}}\,du.
\end{equation}
Indeed, for fixed $\eta>0$ the vertical sides tend to zero as their
real parts tend to infinity.

Apply \eqref{eq:shifted-K} with $x=2\pi py$ and $x=2\pi qy$ and insert the
result into the definition of $\cB_{p,q}(t)$.  Since
$\supp\psi\subset[y_0,y_1]\Subset(0,\infty)$, the resulting triple integral
is absolutely convergent for every fixed $\eta>0$.  Thus Fubini's theorem
and the change of variables $\xi=u+v$ give
\begin{equation}\label{eq:B-eta-double}
 \cB_{p,q}(t)
 =P_\eta(t)\int_{\mathbb R}e^{it\xi}
 \mathscr A^{(\eta)}_{p,q}(\xi)\,d\xi,
 \qquad
 P_\eta(t):=4e^{-(\pi-2\eta)t}\cosh(\pi t),
\end{equation}
where
\begin{align*}
 \mathscr A^{(\eta)}_{p,q}(\xi)
 :=\int_{\mathbb R}\int_0^\infty &\psi(y)
 e^{-2\pi y\sin\eta\{p\cosh u+q\cosh(\xi-u)\}}\\
 &\times e^{-2\pi iy\cos\eta\{p\sinh u+q\sinh(\xi-u)\}}\,dy\,du.
\end{align*}

It remains to let $\eta\downarrow0$.  Take
$\phi\in C_c^\infty(0,\infty)$ and put
\[
 \Phi(s)=\int_0^\infty\phi(t)e^{its}\,dt.
\]
Then $\Phi$ is Schwartz.  Multiplying \eqref{eq:B-eta-double} by
$\phi(t)P_\eta(t)^{-1}$ and integrating in $t$ gives
\begin{equation}\label{eq:test-function-eta}
 \int_0^\infty\phi(t)P_\eta(t)^{-1}\cB_{p,q}(t)\,dt
 =\int_{\mathbb R}\Phi(\xi)\mathscr A^{(\eta)}_{p,q}(\xi)\,d\xi.
\end{equation}

Write
\[
 F_\xi(u)=p\sinh u+q\sinh(\xi-u),
\]
and put
\[
 A_\eta=2\pi\sin\eta\{p\cosh u+q\cosh(\xi-u)\},
 \qquad
 B_\eta=2\pi\cos\eta\,F_\xi(u).
\]
Since $\supp\psi\subset[y_0,y_1]\Subset(0,\infty)$, integration by
parts $M$ times in $y$ gives, for every $M\ge0$,
\[
 \left|\int_0^\infty\psi(y)e^{-A_\eta y}e^{-iB_\eta y}\,dy\right|
 \ll_{M,\psi}\langle B_\eta\rangle^{-M}.
\]
Indeed, differentiating $\psi(y)e^{-A_\eta y}$ produces terms bounded
by $A_\eta^k e^{-A_\eta y_0}$, and
$A_\eta^k e^{-A_\eta y_0}\ll_{k,y_0}1$.  Hence, uniformly for
$0<\eta\le\eta_0<\pi/2$,
\begin{equation}\label{eq:Psi-uniform-majorant}
 \left|\int_0^\infty
 \psi(y)e^{-2\pi y\sin\eta\{p\cosh u+q\cosh(\xi-u)\}}
 e^{-2\pi i y\cos\eta F_\xi(u)}\,dy\right|
 \ll_{M,\psi,\eta_0}\langle F_\xi(u)\rangle^{-M}
 \qquad(M\ge0),
\end{equation}
since $\cos\eta\ge\cos\eta_0>0$.

The constants in these estimates depend only on $M$, $\psi$, and $\eta_0$,
so the following majorant is a function of a single real variable.
Choose constants $C_M$ for which the right-hand side of
\eqref{eq:Psi-uniform-majorant} is bounded by
$C_M\langle F_\xi(u)\rangle^{-M}$, and put
\[
 g(s):=\inf_{M\ge0} C_M\langle s\rangle^{-M}.
\]
Then $g$ is upper semicontinuous, $g(s)\ll_M\langle s\rangle^{-M}$ for every $M$, and
\eqref{eq:Psi-uniform-majorant} is bounded by $g(F_\xi(u))$.
By Lemma~\ref{lem:schwartz-majorant}, there is a nonnegative
$\rho\in\mathcal S(\mathbb R)$, depending only on $\psi$ and $\eta_0$,
such that $g\le\rho$.  Thus
\begin{equation}\label{eq:Psi-Schwartz-majorant}
 \left|\int_0^\infty
 \psi(y)e^{-2\pi y\sin\eta\{p\cosh u+q\cosh(\xi-u)\}}
 e^{-2\pi i y\cos\eta F_\xi(u)}\,dy\right|
 \le \rho(F_\xi(u)).
\end{equation}
By Lemmas~\ref{lem:A-exact} and \ref{lem:model-integrals}, applied with
$\varphi=\rho$,
\[
 \xi\longmapsto \int_{\mathbb R}\rho(F_\xi(u))\,du
\]
is locally integrable, with at most logarithmic singularities at
$\xi=\pm\vartheta$, and has rapid decay as $|\xi|\to\infty$.
Since $\Phi$ is Schwartz,
\[
 \int_{\mathbb R}|\Phi(\xi)|
 \left(\int_{\mathbb R}\rho(F_\xi(u))\,du\right)d\xi<\infty.
\]
Thus \eqref{eq:Psi-Schwartz-majorant} gives an $\eta$-independent
integrable majorant in \eqref{eq:test-function-eta}.  Dominated convergence yields
\begin{align}
 \lim_{\eta\downarrow0}
 \int_{\mathbb R}\Phi(\xi)\mathscr A^{(\eta)}_{p,q}(\xi)\,d\xi
 &=\int_{\mathbb R}\Phi(\xi)
   \left(\int_{\mathbb R}\widehat\psi(F_\xi(u))\,du\right)d\xi\notag\\
 &=\int_{\mathbb R}\Phi(\xi)\mathscr A_{p,q}(\xi)\,d\xi.
                                                        \label{eq:test-function-limit}
\end{align}
Since $P_\eta^{-1}\to P_0^{-1}$ uniformly on $\supp\phi$, with
$P_0(t)=2(1+e^{-2\pi t})$, equations \eqref{eq:test-function-eta} and
\eqref{eq:test-function-limit} give
\[
 \frac{\cB_{p,q}(t)}{P_0(t)}
 =\int_{\mathbb R}e^{it\xi}\mathscr A_{p,q}(\xi)\,d\xi
\]
in the distributional sense on $(0,\infty)$.  By
\eqref{eq:A-weighted-L1} with $m=0$, the Fourier integral on the right is
continuous; the left-hand side is continuous as well.  Thus the identity
holds pointwise, proving \eqref{eq:B-Fourier-exact}.

Finally, \eqref{eq:A-weighted-L1} permits differentiation under the integral:
\[
 \frac{d^m}{dt^m}
 \int_{\mathbb R}e^{it\xi}\mathscr A_{p,q}(\xi)\,d\xi
 =\int_{\mathbb R}(i\xi)^me^{it\xi}
 \mathscr A_{p,q}(\xi)\,d\xi.
\]
Applying Leibniz' rule to \eqref{eq:B-Fourier-exact} gives
\eqref{eq:B-derivative-exact}.
\end{proof}

\subsection{Estimates near and away from the zeros of \texorpdfstring{\(\mathcal D_{p,q}\)}{D(p,q)}}

\begin{lem}\label{lem:BV-Fourier}
Let $I=[\alpha,\beta]$, $H\in BV(I)$, and $m\in C^1(I)$.  Then
\begin{equation}\label{eq:BV-product}
 \operatorname{Var}_{I}(mH)
 \le \|m\|_{L^\infty(I)}\operatorname{Var}_{I}(H)
 +\|H\|_{L^\infty(I)}\int_I|m'(s)|\,ds.
\end{equation}
For every $t\ne0$,
\begin{align}
 \left|\int_\alpha^\beta e^{its}m(s)H(s)\,ds\right|
 &\le \frac1{|t|}\left(
  \{2\|m\|_\infty+\|m'\|_{L^1(I)}\}\|H\|_\infty
  +\|m\|_\infty\operatorname{Var}_{I}(H)\right).       \label{eq:BV-Fourier}
\end{align}
\end{lem}

\begin{proof}
For a partition $\alpha=x_0<\cdots<x_n=\beta$, use
\[
 |m(x_k)H(x_k)-m(x_{k-1})H(x_{k-1})|
 \le |m(x_k)|\,|H(x_k)-H(x_{k-1})|
 +|H(x_{k-1})|\,|m(x_k)-m(x_{k-1})|.
\]
After summing and taking the supremum, the identity
$\operatorname{Var}_{I}(m)=\int_I|m'|$ proves
\eqref{eq:BV-product}.

Let $F=mH$.  Its distributional derivative is a finite complex measure
$dF$ with $|dF|(I)=\operatorname{Var}_{I}(F)$.  Stieltjes integration by
parts gives
\[
 \int_\alpha^\beta e^{its}F(s)\,ds
 =\frac{e^{it\beta}F(\beta-)-e^{it\alpha}F(\alpha+)}{it}
 -\frac1{it}\int_{(\alpha,\beta]}e^{its}\,dF(s).
\]
Taking absolute values and using $|e^{its}|=1$ yields
\[
 \left|\int_\alpha^\beta e^{its}F(s)\,ds\right|
 \le\frac{2\|F\|_\infty+\operatorname{Var}_{I}(F)}{|t|}.
\]
Together with $\|F\|_\infty\le\|m\|_\infty\|H\|_\infty$ and
\eqref{eq:BV-product}, this proves \eqref{eq:BV-Fourier}.
\end{proof}

\begin{lem}
\label{lem:xi-integral}
Fix $0<\eps<1$ and $T\ge2$, and assume
\[
 1\le p<q\le C_\psi T,
 \qquad
 1\le q^2-p^2\ll T^{1+\eps},
\]
and put
\begin{equation*}
 \vartheta=\log(q/p),\qquad \Delta=T^{-1+\eps}.
\end{equation*}
For every fixed integer $r\ge0$ and $t\asymp T$,
\begin{equation}\label{eq:xi-integral}
 \int_{\mathbb R}(i\xi)^re^{it\xi}\mathscr A_{p,q}(\xi)\,d\xi
 \ll_{r,\psi,\eps}\frac{\log T}{T}
 (\vartheta+\Delta)^r.
\end{equation}
\end{lem}

\begin{proof}
Choose $\chi\in C_c^\infty((-2,2))$ with $\chi=1$ on $[-1,1]$, and set
\[
 \chi_+(\xi)=\chi((\xi-\vartheta)/\Delta),\qquad
 \chi_-(\xi)=\chi((\xi+\vartheta)/\Delta).
\]
Use the partition
\begin{equation*}
 \eta_+=\chi_+,\qquad
 \eta_-=(1-\chi_+)\chi_-,\qquad
 \eta_0=(1-\chi_+)(1-\chi_-),
 \qquad \eta_++\eta_-+\eta_0=1.
\end{equation*}
For every $\ell\ge0$,
\begin{equation}\label{eq:cutoff-derivatives}
 |\partial_\xi^\ell\eta_\pm(\xi)|
 +|\partial_\xi^\ell\eta_0(\xi)|
 \ll_\ell\Delta^{-\ell}.
\end{equation}
Since $t\asymp T$,
\[
 t\Delta\asymp T^{\eps}\ge1
\]
for large $T$; bounded $T$ is absorbed by increasing the implicit constant.

We first estimate the two terms supported near the zeros $\pm\vartheta$.
There we use the local logarithmic representation \eqref{eq:A-local-log}:
the singular logarithm is treated directly, while the remaining term is
handled by the $BV$ estimate \eqref{eq:H-BV}.  Suppose first that the two
$\Delta$-neighbourhoods are separated, so that
\[
 \vartheta\ge4\Delta,\qquad
 \xi=\vartheta+s,\qquad |s|\le2\Delta.
\]
Using
\[
 \mathcal D_{p,q}(\xi)
 =
 4pq\sinh\frac{\vartheta+\xi}{2}
     \sinh\frac{\vartheta-\xi}{2},
 \qquad
 q^2-p^2=2pq\sinh\vartheta,
\]
we obtain
\begin{equation*}
 |\mathcal D_{p,q}(\vartheta+s)|
 =(q^2-p^2)|s|e^{b_+(s)},
\end{equation*}
where
\begin{equation*}
 b_+(s)
 =
 \log\!\left(\frac{2|\sinh(s/2)|}{|s|}\right)
 +
 \log\!\left(
  \frac{\sinh(\vartheta+s/2)}{\sinh\vartheta}
 \right),
\end{equation*}
with $b_+(0)=0$.  Thus
\begin{equation*}
 \log|\mathcal D_{p,q}(\vartheta+s)|
 =
 \log(q^2-p^2)+\log|s|+b_+(s).
\end{equation*}

Write
\[
 b_{+,1}(s)=
 \log\!\left(\frac{2|\sinh(s/2)|}{|s|}\right),
 \qquad
 b_{+,2}(s)=
 \log\!\left(
  \frac{\sinh(\vartheta+s/2)}{\sinh\vartheta}
 \right).
\]
On $|s|\le1/2$,
\[
 b_{+,1}(s)=O(s^2),\qquad b_{+,1}'(s)=O(|s|).
\]
Moreover,
\[
 b_{+,2}'(s)=\frac12\coth(\vartheta+s/2).
\]
Since
\[
 \vartheta+s/2\ge\vartheta-\Delta
 \ge\frac34\vartheta,
\]
if $0<\vartheta\le1$, then
\[
 |b_{+,2}'(s)|\ll\vartheta^{-1},
 \qquad
 |b_{+,2}(s)|
 \ll\frac{\Delta}{\vartheta}\ll1.
\]
If $\vartheta>1$, the hyperbolic cotangent is bounded, and the same
conclusion follows.  Therefore
\begin{equation}\label{eq:bplus-bound}
 \|b_+\|_{L^\infty(|s|\le2\Delta)}
 +
 \Delta
 \|b_+'\|_{L^\infty(|s|\le2\Delta)}
 \ll1.
\end{equation}
At the zero $-\vartheta$, write $s=\xi+\vartheta$.  Evenness of
$\mathcal D_{p,q}$ gives the same formula and estimate, with
$b_-(s)=b_+(-s)$.

We now estimate the contribution of these two neighbourhoods.  On the
support of either \(\eta_+\) or \(\eta_-\),
\begin{equation}\label{eq:w-power-local}
 |\xi|^r\le C_r(\vartheta+\Delta)^r.
\end{equation}
If
\[
 m_{r,+}(s)=(i(\vartheta+s))^r\eta_+(\vartheta+s),
\]
then
\begin{equation}\label{eq:threshold-multiplier}
 \|m_{r,+}\|_\infty+
 \Delta\|m_{r,+}'\|_\infty
 \ll_r(\vartheta+\Delta)^r.
\end{equation}
The same estimate holds near the zero $-\vartheta$.

Let $c\in C^1([-2\Delta,2\Delta])$.
On $|s|\le t^{-1}$,
\begin{equation*}
 \int_{|s|\le t^{-1}}|\log|s||\,ds
 =2t^{-1}(1+\log t)\ll t^{-1}\log(2+t).
\end{equation*}
On the two intervals $t^{-1}\le|s|\le2\Delta$, integration
by parts yields
\[
 \int e^{its}\log|s|c(s)\,ds
 =
 \frac{e^{its}\log|s|c(s)}{it}
 -
 \frac1{it}\int e^{its}
 \left(\frac{c(s)}s+\log|s|c'(s)\right)\,ds.
\]
Hence
\begin{equation}\label{eq:log-ibp}
 \left|\int_{|s|\le2\Delta}
 e^{its}\log|s|c(s)\,ds\right|
 \ll
 t^{-1}\log(2+t)
 \left\{\|c\|_\infty+
 \Delta\|c'\|_\infty\right\}.
\end{equation}

Apply \eqref{eq:A-local-log}.  At the zero $\vartheta$ put
\[
 \mathcal H_+(s)=H_+(\vartheta+s),
 \qquad |s|\le2\Delta.
\]
Translation invariance of total variation, \eqref{eq:H-BV},
\eqref{eq:threshold-multiplier}, and Lemma~\ref{lem:BV-Fourier} give
\begin{align*}
 &\left|\int_{|s|\le2\Delta}
 e^{its}m_{r,+}(s)\mathcal H_+(s)\,ds\right| \\
 &\quad\ll
 \frac1t
 \{\|m_{r,+}\|_\infty+
   \Delta\|m_{r,+}'\|_\infty\}
 \{\|H_+\|_\infty+\operatorname{Var}_{I_+}(H_+)\} \\
 &\quad\ll_{r,\psi,\eps}
 T^{-1}(\log T)(\vartheta+\Delta)^r.
\end{align*}
The term $\log|s|$ is controlled by \eqref{eq:log-ibp} and
\eqref{eq:threshold-multiplier}.

Since $m_{r,+}(\pm2\Delta)=0$, one ordinary integration by
parts gives
\begin{align*}
 &\left|
 \int_{|s|\le2\Delta}
 e^{its}m_{r,+}(s)
 \{\log(q^2-p^2)+b_+(s)\}\,ds
 \right| \\
 &\quad\le
 \frac1t\int_{|s|\le2\Delta}
 \left|
 \frac{d}{ds}
 \left[
 m_{r,+}(s)\{\log(q^2-p^2)+b_+(s)\}
 \right]\right|\,ds \\
 &\quad\ll
 \frac1t\left[
 \Delta\|m_{r,+}'\|_\infty
 \{\log T+\|b_+\|_\infty\}
 +
 \Delta\|m_{r,+}\|_\infty
 \|b_+'\|_\infty
 \right] \\
 &\quad\ll_{r,\psi,\eps}
 T^{-1}(\log T)
 (\vartheta+\Delta)^r.
\end{align*}
Here
\[
 1\le q^2-p^2\ll T^{1+\eps},
 \qquad
 |\log(q^2-p^2)|\ll_\eps\log T,
\]
and \eqref{eq:bplus-bound} were used.  The argument at the zero $-\vartheta$ is identical.

It remains to treat the case in which the two $\Delta$-neighbourhoods
overlap.  If $\vartheta<4\Delta$, the second zero is not absorbed into
$b_+$.  Instead,
\[
 \log|\mathcal D_{p,q}(\xi)|
 =
 \log(4pq)
 +
 \log\left|\sinh\frac{\xi+\vartheta}{2}\right|
 +
 \log\left|\sinh\frac{\xi-\vartheta}{2}\right|.
\]
Both zeros lie in $|\xi|\le4\Delta$, and the support of
$\eta_++\eta_-$ lies in $|\xi|\le6\Delta$.  Write
\[
 \log\left|\sinh\frac{\xi-\xi_0}{2}\right|
 =
 \log|\xi-\xi_0|
 +
 \log\left|
 \frac{\sinh((\xi-\xi_0)/2)}{\xi-\xi_0}
 \right|,
 \qquad \xi_0=\pm\vartheta.
\]
The second term and its first derivative are uniformly bounded, while the
first is handled by \eqref{eq:log-ibp}.  Together with
\eqref{eq:w-power-local} and the BV estimate, this gives
\[
 T^{-1}(\log T)(\vartheta+\Delta)^r.
\]

Thus
\begin{equation}\label{eq:threshold-piece-final}
 \left|\int e^{it\xi}(i\xi)^r
 (\eta_++\eta_-)\mathscr A_{p,q}(\xi)\,d\xi\right|
 \ll_{r,\psi,\eps}
 T^{-1}(\log T)(\vartheta+\Delta)^r.
\end{equation}

We now turn to the complementary term containing $\eta_0$.  On its
support we are a distance $\gg\Delta$ from both zeros, so the weighted
$L^1$ derivative estimate \eqref{eq:A-derivative-away} can be combined
with repeated integration by parts.  Leibniz' formula gives
\[
 \partial_\xi^J\{(i\xi)^r\eta_0\mathscr A_{p,q}\}
 =
 \sum_{j_1+j_2+j_3=J}
 c_{j_1,j_2,j_3}
 \partial_\xi^{j_1}(i\xi)^r\,
 \partial_\xi^{j_2}\eta_0\,
 \partial_\xi^{j_3}\mathscr A_{p,q}.
\]
Using \eqref{eq:cutoff-derivatives} and
\eqref{eq:A-derivative-away},
\begin{equation*}
 \int_{\mathbb R}
 \left|\partial_\xi^J\{(i\xi)^r\eta_0(\xi)
 \mathscr A_{p,q}(\xi)\}\right|\,d\xi
 \ll_{J,r,\psi}
 (\log T)^{r+1}\Delta^{-J}.
\end{equation*}
Here one uses the weighted $L^1$ form of
\eqref{eq:A-derivative-away} with weight at most $(1+|\xi|)^r$; hence the
polynomial factors produced by $\partial_\xi^{j_1}(i\xi)^r$ are absorbed
at the cost of $(\log T)^{r+1}$, while the derivatives falling on
$\eta_0$ and $\mathscr A_{p,q}$ contribute together at most
$\Delta^{-J}$.  After $J$ integrations by parts,
\begin{equation*}
 \left|\int e^{it\xi}(i\xi)^r\eta_0(\xi)
 \mathscr A_{p,q}(\xi)\,d\xi\right|
 \ll_{J,r,\psi}(\log T)^{r+1}(t\Delta)^{-J}
 \ll_{J,r,\psi}(\log T)^{r+1}T^{-J\eps}.
\end{equation*}
Since
\[
 \vartheta+\Delta\ge\Delta=T^{-1+\eps},
\]
we may choose $J=J(r,\eps)$ sufficiently large that
\[
 (\log T)^{r+1}T^{-J\eps}
 \ll T^{-1}(\log T)(\vartheta+\Delta)^r.
\]
Combining this with \eqref{eq:threshold-piece-final} proves
\eqref{eq:xi-integral}.
\end{proof}

\begin{lem}\label{lem:B-derivatives}
Fix $0<\eps<1$. Let $t\asymp T$, put $p=da$, $q=d(a+b)$ and $\rho=b/a$.
Assume the truncation conditions of Lemma~\ref{lem:xi-integral}.
Then, for every fixed integer $r\ge0$,
\begin{equation}\label{eq:derivative-bound}
 \frac{d^r}{dt^r}\cB_{da,d(a+b)}(t)
 \ll_{r,\psi,\eps}\frac{\log T}{T}\Lambda(\rho,T)^r,
\end{equation}
where
\begin{equation}\label{eq:Lambda}
 \Lambda(\rho,T)=T^{-1+\eps}+
 \begin{cases}
  \rho,&0<\rho\le1,\\
  \log(2+\rho),&\rho>1.
 \end{cases}
\end{equation}
\end{lem}

\begin{proof}
Fix the $\eps>0$ in the conclusion and use a smaller parameter
$\eps_0>0$ in the truncation of the preceding section, with
$\eps_0<\eps/2$ at the end.  The original Fourier-index truncation gives
\[
 p=da\le C_\psi T,\qquad q=d(a+b)\le C_\psi T,
\]
and
\begin{equation*}
 1\le q^2-p^2=d^2b(2a+b)\ll T^{1+\eps_0}.
\end{equation*}
Hence Lemma~\ref{lem:xi-integral} applies.  Put
\[
 \vartheta=\log(q/p),\qquad
 \Delta=T^{-1+\eps_0},
\]
and for $0\le m\le r$ define
\[
 \mathcal I_m(t;p,q)
 =
 \int_{\mathbb R}(i\xi)^me^{it\xi}
 \mathscr A_{p,q}(\xi)\,d\xi.
\]
By \eqref{eq:xi-integral},
\begin{equation*}
 \mathcal I_m(t;p,q)
 \ll_{m,\psi,\eps_0}
 T^{-1}(\log T)
 (\vartheta+\Delta)^m.
\end{equation*}

The differentiation formula
\eqref{eq:B-derivative-exact} gives
\begin{align}
 |\partial_t^r\cB_{p,q}(t)|
 &\le
 2(1+e^{-2\pi t})|\mathcal I_r(t;p,q)| \notag\\
 &\quad+
 2e^{-2\pi t}\sum_{\ell=1}^r
 \binom r\ell(2\pi)^\ell
 |\mathcal I_{r-\ell}(t;p,q)|.
                                                        \label{eq:B-derivative-first}
\end{align}
The first term is
\[
 \ll_{r,\psi,\eps_0}
 T^{-1}(\log T)
 (\vartheta+\Delta)^r.
\]
Since
\[
 0<\vartheta=\log(q/p)\le\log(C_\psi T)
 \ll_\psi\log T,
\]
for every fixed $B>0$ the second line in
\eqref{eq:B-derivative-first} is
\begin{align}
 &\ll_{r,\psi,\eps_0}
 e^{-cT}T^{-1}(\log T)^{r+1} \notag\\
 &\ll_{B,r,\psi}T^{-B}.                    \label{eq:B-exp-tail}
\end{align}
Thus
\begin{equation*}
 |\partial_t^r\cB_{p,q}(t)|
 \ll_{r,\psi,\eps_0}
 T^{-1}(\log T)
 (\vartheta+T^{-1+\eps_0})^r
 +O_{B,r,\psi}(T^{-B}).
\end{equation*}

Finally,
\[
 \frac qp=1+\rho,\qquad
 \vartheta=\log(1+\rho).
\]
If $0<\rho\le1$, then directly
\[
 \log(1+\rho)\le\rho,
\]
whereas if $\rho>1$,
\[
 \log(1+\rho)\le\log(2+\rho).
\]
Therefore
\begin{equation*}
 \vartheta+T^{-1+\eps}\le\Lambda(\rho,T).
\end{equation*}
Since $\Lambda(\rho,T)\ge T^{-1+\eps}$, choosing $B$ sufficiently large absorbs
\eqref{eq:B-exp-tail} into
$T^{-1}(\log T)\Lambda(\rho,T)^r$.  Taking $\eps_0<\eps$ proves
\eqref{eq:derivative-bound}.

\end{proof}

\section{The local spectral large sieve with varying coefficients}
\label{sec:sieve}

To apply the spectral large sieve in Lemma~\ref{lem:spectral-sieve} to
the $t$-dependent coefficients, we use the following Fourier-series lemma.

\begin{lem}\label{lem:vector-fourier}
Let $\mathcal N$ be a finite index set, let $J$ be an interval of
length $L\asymp H$ with midpoint $t_J$, and let
\[
 F:J\longrightarrow\ell^2(\mathcal N),\qquad
 F(t)=(F_\nu(t))_{\nu\in\mathcal N},
\]
be smooth and vanish near the endpoints.  Suppose that, for every fixed $r\ge0$,
\begin{equation}\label{eq:vector-derivative}
 \sup_{t\in J}\sum_{\nu\in\mathcal N}|\partial_t^rF_\nu(t)|^2
 \ll_r E\Lambda^{2r}.
\end{equation}
Then $F$ has the Fourier expansion
\begin{equation}\label{eq:vector-Fourier-series}
 F_\nu(t)=\sum_{k\in\mathbb Z}c_{\nu,k}e\!\left(\frac{k(t-t_J)}L\right),
 \qquad
 \sum_{k\in\mathbb Z}\sum_{\nu\in\mathcal N}|c_{\nu,k}|^2\ll E.
\end{equation}
Moreover, for every integer $R\ge1$ and $K\ge1$,
\begin{align}
 \sum_{|k|>K}\sum_{\nu\in\mathcal N}|c_{\nu,k}|^2
 &\ll_R E K^{-2R}(1+H\Lambda)^{2R},                 \label{eq:Fourier-L2-tail}\\
 \sup_{t\in J}\sum_{\nu\in\mathcal N}\left|
  \sum_{|k|>K}c_{\nu,k}e\!\left(\frac{k(t-t_J)}L\right)
 \right|^2
 &\ll_R E K^{1-2R}(1+H\Lambda)^{2R}.               \label{eq:Fourier-sup-tail}
\end{align}
\end{lem}

\begin{proof}
Write $\mathcal H=\ell^2(\mathcal N)$.
Since $\mathcal N$ is finite, $\mathcal H$ is a finite-dimensional
Hilbert space; equivalently, all identities below may be proved componentwise
and then summed over $\nu\in\mathcal N$.  Since $F$ vanishes in
neighbourhoods of both endpoints of $J$, its $L$-periodic extension belongs to
$C^\infty(\mathbb R/L\mathbb Z;\mathcal H)$ and
is smooth across the periodic seam.

\smallskip
\noindent\emph{Step 1: Fourier coefficients and exact expansion.}
For $k\in\mathbb Z$ define the $\mathcal H$-valued coefficient
\begin{equation}\label{eq:vector-Fourier-coeff}
 c_k=(c_{\nu,k})_{\nu\in\mathcal N}
 :=\frac1L\int_{J}F(t)
 e\!\left(-\frac{k(t-t_J)}L\right)\,dt .
\end{equation}
Since the periodic extension of $F$ is smooth, its Fourier series and all
its derivatives converge uniformly.  Hence, on $J$,
\[
 F(t)=\sum_{k\in\mathbb Z}c_k
 e\!\left(\frac{k(t-t_J)}L\right)
\]
as an identity in $\mathcal H$.  Taking the $\nu$-th coordinate gives the first
formula in \eqref{eq:vector-Fourier-series}.

\smallskip
\noindent\emph{Step 2: zeroth-order Parseval identity.}
Parseval on the circle of length $L$, with normalised Haar measure $dt/L$,
gives
\begin{align*}
 \sum_{k\in\mathbb Z}\|c_k\|_{\mathcal H}^2
 &=\frac1L\int_{J}\|F(t)\|_{\mathcal H}^2\,dt \\
 &=\frac1L\int_{J}\sum_{\nu\in\mathcal N}|F_\nu(t)|^2\,dt
 \ll E.
\end{align*}
The last bound is \eqref{eq:vector-derivative} with $r=0$.  Since
$\|c_k\|_{\mathcal H}^2=\sum_{\nu\in\mathcal N}|c_{\nu,k}|^2$, this is the second assertion in
\eqref{eq:vector-Fourier-series}.

\smallskip
\noindent\emph{Step 3: weighted Fourier energy.}
Integrate \eqref{eq:vector-Fourier-coeff} by parts $R$ times.  Every boundary
term vanishes because $F$ is identically zero near both endpoints.  Thus, for
every $k\in\mathbb Z$,
\begin{equation*}
 \left(\frac{2\pi i k}{L}\right)^Rc_k
 =\frac1L\int_{J}\partial_t^RF(t)
 e\!\left(-\frac{k(t-t_J)}L\right)\,dt .
\end{equation*}
Equivalently, the $k$-th Fourier coefficient of $\partial_t^RF$ is
$(2\pi i k/L)^Rc_k$.  A second application of vector-valued Parseval yields
\begin{align*}
 \sum_{k\in\mathbb Z}
 \left(\frac{2\pi|k|}{L}\right)^{2R}\|c_k\|_{\mathcal H}^2
 &=\frac1L\int_{J}
   \|\partial_t^RF(t)\|_{\mathcal H}^2\,dt \\
 &\ll_R E\Lambda^{2R}.
\end{align*}
Consequently,
\begin{equation}\label{eq:weighted-Fourier-energy}
 \sum_{k\in\mathbb Z}|k|^{2R}\|c_k\|_{\mathcal H}^2
 \ll_R \left(\frac{L}{2\pi}\right)^{2R}E\Lambda^{2R}
 \ll_R E(1+H\Lambda)^{2R},
\end{equation}
because $L\asymp H$.  The factor $1+H\Lambda$ also covers the case $H\Lambda<1$.  In the application in Section~\ref{sec:offdiag}, multiplication
by a cutoff supported in \(J\) replaces \(\Lambda\) by
\(\Lambda+H^{-1}\), for which the same expression is exactly of size
$1+H\Lambda$.

\smallskip
\noindent\emph{Step 4: the $L^2$ Fourier tail.}
For $|k|>K$ one has $1\le K^{-2R}|k|^{2R}$.  Hence
\begin{align*}
 \sum_{|k|>K}\sum_{\nu\in\mathcal N}|c_{\nu,k}|^2
 &=\sum_{|k|>K}\|c_k\|_{\mathcal H}^2\\
 &\le K^{-2R}\sum_{k\in\mathbb Z}
       |k|^{2R}\|c_k\|_{\mathcal H}^2\\
 &\ll_R EK^{-2R}(1+H\Lambda)^{2R}
\end{align*}
by \eqref{eq:weighted-Fourier-energy}.  This proves
\eqref{eq:Fourier-L2-tail}.

\smallskip
\noindent\emph{Step 5: a tail uniform in the spectral parameter.}
Set
\[
 \mathcal R_K(t)
 :=\sum_{|k|>K}c_k
 e\!\left(\frac{k(t-t_J)}L\right)\in\mathcal H.
\]
For each fixed $t\in J$, the triangle inequality followed by Cauchy's
inequality for the sequences $|k|^{-R}$ and
$|k|^R\|c_k\|_{\mathcal H}$ gives
\begin{align}
 \|\mathcal R_K(t)\|_{\mathcal H}^2
 &\le
 \left(\sum_{|k|>K}|k|^{-R}
       |k|^R\|c_k\|_{\mathcal H}\right)^2 \notag\\
 &\le
 \left(\sum_{|k|>K}|k|^{-2R}\right)
 \left(\sum_{|k|>K}|k|^{2R}\|c_k\|_{\mathcal H}^2\right).
                                                               \label{eq:uniform-tail-cauchy}
\end{align}
For $R\ge1$,
\[
 \sum_{|k|>K}|k|^{-2R}\ll_R K^{1-2R}.
\]
Combining this estimate with \eqref{eq:weighted-Fourier-energy} in
\eqref{eq:uniform-tail-cauchy}, we obtain
\[
 \sup_{t\in J}\|\mathcal R_K(t)\|_{\mathcal H}^2
 \ll_R EK^{1-2R}(1+H\Lambda)^{2R}.
\]
Finally,
\[
 \|\mathcal R_K(t)\|_{\mathcal H}^2
 =\sum_{\nu\in\mathcal N}\left|
   \sum_{|k|>K}c_{\nu,k}
   e\!\left(\frac{k(t-t_J)}L\right)\right|^2,
\]
which is \eqref{eq:Fourier-sup-tail}.
\end{proof}

\section{Proof of the main theorem}\label{sec:offdiag}

We first prove the following weighted second-moment estimate for the
off-diagonal term.
\begin{theorem}\label{thm:offdiag}
Let $S_j(\psi)$ be defined by \eqref{eq:S-hecke}.  For $T\ge2$,
$1\le U\le T$, and any $\eps>0$,
\begin{equation}\label{eq:Svariance}
 \sum_{\substack{j\ {\rm even}\\T<t_j\le T+U}}
 L_j|S_j(\psi)|^2
 \ll_{\psi,\eps}T^{1+\eps}.
\end{equation}
\end{theorem}

\begin{proof}
Recall \eqref{eq:S-hecke}:
\begin{equation*}
 S_j(\psi)=\frac1{L_j}\sum_{d,a,b\ge1}
 \lambda_j(a(a+b))\cB_{da,d(a+b)}(t_j).
\end{equation*}

First, we divide the threefold sum into dyadic blocks and estimate the
coefficients of the resulting Hecke polynomials.
We decompose the retained triples,
subject to the Fourier-index and shift
truncations from Section~\ref{sec:trunc}, into dyadic ranges
\[
 d\asymp D,\qquad a\asymp X,\qquad b\asymp Y.
\]
Call the resulting block $\mathfrak B=(D,X,Y)$.  In all sums over this
block we retain only triples satisfying
\[
 da,\ d(a+b)\le C_\psi T,\qquad
 d^2b(2a+b)\ll T^{1+\beta}.
\]
For each Hecke index $A$, let $\mathscr T_A(\mathfrak B)$ denote the
triples in this block that satisfy $a(a+b)=A$ and the truncation
conditions above. Define
\begin{equation*}
 C_A^{\mathfrak B}(t)=
 \sum_{(d,a,b)\in\mathscr T_A(\mathfrak B)}
 \cB_{da,d(a+b)}(t).
\end{equation*}
The Hecke index has length
\begin{equation*}
 A\asymp N:=X(X+Y).
\end{equation*}

For every fixed integer $r\ge0$, we claim that
\begin{equation}\label{eq:energy}
 \sum_A\left|\partial_t^rC_A^{\mathfrak B}(t)\right|^2
 \ll_{r,\eps}
 E_{\mathfrak B}\Lambda(Y/X,T)^{2r},
 \qquad
 E_{\mathfrak B}=T^\eps\frac{D^2XY}{T^2}.
\end{equation}

We consider only nonempty blocks left after the truncation in
Section~\ref{sec:trunc}.  To keep the bookkeeping transparent, apply
Lemma~\ref{lem:B-derivatives} and the divisor bound with $\eps/10$ in place of
$\eps$; the finitely many resulting powers $T^{\eps/10}$ are absorbed into
$T^\eps$ at the end.

Then
\begin{equation}\label{eq:block-energy-expand}
 \partial_t^rC_A^{\mathfrak B}(t)
 =\sum_{(d,a,b)\in\mathscr T_A(\mathfrak B)}
   \partial_t^r\mathcal B_{da,d(a+b)}(t).
\end{equation}

\par
If $(d,a,b)\in\mathscr T_A(\mathfrak B)$, then
\[
 a\mid A,\qquad b=A/a-a.
\]
Thus $b$ is uniquely determined once the divisor $a$ is chosen.  There are
$O(D)$ integral choices for $d$, and hence
\begin{equation*}
 \#\mathscr T_A(\mathfrak B)\ll D\,\tau(A).
\end{equation*}
The truncation in Lemma~\ref{lem:K-tail} gives
$da,d(a+b)\ll_\psi T$ for every retained triple, and therefore
$A=a(a+b)\ll_\psi T^2$.  The classical divisor estimate
$\tau(A)\ll_\eps A^{\eps/20}$ now yields
\begin{equation}\label{eq:block-multiplicity-final}
 \#\mathscr T_A(\mathfrak B)\ll_{\psi,\eps}D\,T^{\eps/10}.
\end{equation}

\par
For $a\asymp X$ and $b\asymp Y$,
\[
 \frac b a\asymp\frac YX.
\]
The piecewise definition \eqref{eq:Lambda} implies
\begin{equation*}
 \Lambda(b/a,T)\ll \Lambda(Y/X,T).
\end{equation*}
Indeed, if $Y/X\le1/4$, then $b/a\le1$ and
$b/a\asymp Y/X$; if $Y/X\ge4$, then $b/a>1$ and
$\log(2+b/a)\asymp\log(2+Y/X)$; on the remaining compact interval
$1/4\le Y/X\le4$, both sides are comparable with a positive absolute
constant.  Lemma~\ref{lem:B-derivatives} therefore gives, uniformly throughout
the block,
\begin{equation}\label{eq:block-B-pointwise}
 \left|\partial_t^r\mathcal B_{da,d(a+b)}(t)\right|
 \ll_{r,\psi,\eps}
 T^{-1}(\log T)\Lambda(Y/X,T)^r.
\end{equation}

\par
Equations \eqref{eq:block-energy-expand} and
\eqref{eq:block-multiplicity-final}, together with Cauchy's inequality, give
\begin{align*}
 \left|\partial_t^rC_A^{\mathfrak B}(t)\right|^2
 &\le \#\mathscr T_A(\mathfrak B)
 \sum_{(d,a,b)\in\mathscr T_A(\mathfrak B)}
 \left|\partial_t^r\mathcal B_{da,d(a+b)}(t)\right|^2 \\
 &\ll_{\psi,\eps}DT^{\eps/10}
 \sum_{(d,a,b)\in\mathscr T_A(\mathfrak B)}
 \left|\partial_t^r\mathcal B_{da,d(a+b)}(t)\right|^2 .
\end{align*}
When this is summed over $A$, every triple $(d,a,b)$ occurs exactly once,
namely for $A=a(a+b)$.  Hence
\begin{align*}
 \sum_A\left|\partial_t^rC_A^{\mathfrak B}(t)\right|^2
 &\ll_{\psi,\eps}DT^{\eps/10}
 \sum_{\substack{d\asymp D\\a\asymp X\\b\asymp Y}}
 \left|\partial_t^r\mathcal B_{da,d(a+b)}(t)\right|^2 \\
 &\ll_{r,\psi,\eps}
 DT^{\eps/10}\,(DXY)\,
 T^{-2+\eps/5}\Lambda(Y/X,T)^{2r}.
\end{align*}
The dyadic box contains $O(DXY)$ integral triples, and the second line uses
\eqref{eq:block-B-pointwise}.  Absorbing the finitely many $\eps/10$ losses
into $T^\eps$, we obtain
\[
 \sum_A|\partial_t^rC_A^{\mathfrak B}(t)|^2
 \ll_{r,\psi,\eps}
 T^\eps\frac{D^2XY}{T^2}\Lambda(Y/X,T)^{2r},
\]
which proves \eqref{eq:energy}.

Now the truncation \eqref{eq:trunc-main} becomes
\begin{equation*}
 D^2Y(2X+Y)\ll T^{1+\eps}.
\end{equation*}
Next, we fix a non-empty block $\mathfrak B=(D,X,Y)$ and put
\[
 N=X(X+Y),\qquad \rho=Y/X,\qquad \Lambda=\Lambda(\rho,T).
\]
Define
\begin{equation*}
 \mathcal Q_{\mathfrak B}(T,U)
 =\sum_{\substack{j\ {\rm even}\\T<t_j\le T+U}}\frac1{L_j}
 \left|\sum_{A\asymp N}\lambda_j(A)C_A^{\mathfrak B}(t_j)\right|^2.
\end{equation*}
Let
\begin{equation*}
 H=\begin{cases}
  \min(T^{1-\beta},X/Y),&Y\le X,\\
  1,&Y>X,
 \end{cases}
 \qquad h=\min(U,H).
\end{equation*}
We prove
\begin{equation*}
 \mathcal Q_{\mathfrak B}(T,U)
 \ll T^{O(\beta)}\left(TU+N\max\left(1,\frac UH\right)\right)
 \frac{D^2XY}{T^2}.
\end{equation*}

For this fixed block, choose a smooth square partition of unity
$\{\eta_J\}_{J\in\mathscr J}$ on $[T,T+U]$, indexed by a finite
collection $\mathscr J=\mathscr J(\mathfrak B)$ of short spectral intervals, such that
\begin{equation*}
 \sum_{J\in\mathscr J}\eta_J(t)^2=1,\quad \supp\eta_J\Subset J,\quad
 |J|\asymp h,\quad \eta_J^{(r)}(t)\ll_rh^{-r}.
\end{equation*}
The overlap is bounded and the number of intervals is
\begin{equation*}
 O(U/h)=O\!\left(\max\left(1,\frac UH\right)\right).
\end{equation*}
Put $F_{A,J}(t)=\eta_J(t)C_A^{\mathfrak B}(t)$.  The square partition of unity gives
\begin{align*}
 \mathcal Q_{\mathfrak B}(T,U)
 &=\sum_{J\in\mathscr J}\sum_{\substack{j\ {\rm even}\\T<t_j\le T+U}}\frac1{L_j}
 \left|\sum_A\lambda_j(A)F_{A,J}(t_j)\right|^2\\
 &\le\sum_{J\in\mathscr J}\sum_{j\ {\rm even}}\frac1{L_j}
 \left|\sum_A\lambda_j(A)F_{A,J}(t_j)\right|^2.
\end{align*}
Leibniz' rule and \eqref{eq:energy} give
\begin{equation*}
 \sum_A|\partial_t^rF_{A,J}(t)|^2
 \ll_{r,\psi,\beta}E_{\mathfrak B}(\Lambda+h^{-1})^{2r},
 \qquad E_{\mathfrak B}=T^{O(\beta)}\frac{D^2XY}{T^2}.
\end{equation*}
Since $h\le H$,
\begin{equation*}
 h\Lambda\ll\begin{cases}1,&Y\le X,\\1+\log T,&Y>X.\end{cases}
\end{equation*}
Lemma~\ref{lem:vector-fourier}, with derivative factor \(\Lambda+h^{-1}\),
and $K=T^\beta$ (with the Fourier-tail order chosen sufficiently large depending on $\beta$), gives
\begin{equation*}
 F_{A,J}(t)=\sum_{|k|\le K}c_{A,k,J}
 e\!\left(\frac{k(t-t_J)}{L_J}\right)+R_{A,J}(t),
 \qquad L_J\asymp h,
\end{equation*}
with
\begin{align*}
 \sum_{k,A}|c_{A,k,J}|^2&\ll E_{\mathfrak B},\\
 \sup_{t\in J}\sum_A|R_{A,J}(t)|^2&\ll T^{-200}E_{\mathfrak B}.
\end{align*}
For each fixed $k$, the coefficients $c_{A,k,J}$ are independent of $t_j$, while the exponential factor has absolute value one. Cauchy's inequality in $k$ and Lemma~\ref{lem:spectral-sieve} yield
\begin{equation*}
 \sum_{\substack{j\ {\rm even}\\t_j\in J}}\frac1{L_j}
 \left|\sum_{|k|\le K}e\!\left(\frac{k(t_j-t_J)}{L_J}\right)
 \sum_Ac_{A,k,J}\lambda_j(A)\right|^2
 \ll T^{O(\beta)}(Th+N)E_{\mathfrak B}.
\end{equation*}
The uniform tail, the elementary Hecke bound, and Weyl's upper bound on
$[T,2T]$ make the total tail $O(T^{-190})$.  Summing over these
intervals gives
\begin{align*}
 \mathcal Q_{\mathfrak B}(T,U)
 &\ll T^{O(\beta)}\frac Uh(Th+N)E_{\mathfrak B}+O(T^{-190})\\
 &\ll T^{O(\beta)}\left(TU+N\max\left(1,\frac UH\right)\right)
 \frac{D^2XY}{T^2}.
\end{align*}

Note that the block constraint gives
\begin{equation*}
 TU\frac{D^2XY}{T^2}\ll T^{O(\beta)}U,
\end{equation*}
and
\begin{equation*}
 N\frac{D^2XY}{T^2}\ll T^{O(\beta)}\frac{X^2}{T}.
\end{equation*}
Thus it remains to estimate
\[
 \frac{X^2}{T}\max\left(1,\frac UH\right).
\]
If $Y>X$, then $H=1$ and $X^2<Y^2\ll T^{1+O(\beta)}$, giving
$O(T^{O(\beta)}U)$.  Suppose $Y\le X$ and write $\rho=Y/X$.  If
$T^{-1+\beta}\le\rho\le1$, then $H=\rho^{-1}$ and
$\rho X^2\ll T^{1+O(\beta)}$.  For $U\rho\ge1$ the expression is
$O(T^{O(\beta)}U)$; for $U\rho<1$ it is
$O(T^{O(\beta)}\rho^{-1})=O(T^{1+O(\beta)})$.
If $0<\rho<T^{-1+\beta}$, then $H=T^{1-\beta}$.  Since $d\ge1$ and the
Fourier-index truncation gives $dX\ll_\psi T$, we have $X\ll_\psi T$, and hence
\begin{equation*}
 \frac{X^2}{T}\max\left(1,\frac U{T^{1-\beta}}\right)
 \ll_\psi T^{1+\beta}.
\end{equation*}
In particular, the range $U\rho<1$ gives a bound of order
$T^{1+O(\beta)}$, rather than $T^{O(\beta)}U$.
Since $1\le U\le T$, the preceding cases yield, uniformly for
all non-empty blocks,
\begin{equation*}
 \mathcal Q_{\mathfrak B}(T,U)
 \ll T^{1+O(\beta)}.
\end{equation*}

For each block \(\mathfrak B\), put
\[
 P_{j,\mathfrak B}
 =\sum_A\lambda_j(A)C_A^{\mathfrak B}(t_j).
\]
Then
\[
 S_j(\psi)
 =\frac1{L_j}\sum_{\mathfrak B}P_{j,\mathfrak B}
   +\mathcal E_j^{\mathrm{far}}
  +\mathcal E_j^{\mathrm{shift}}.
\]

There are $O_\psi((\log T)^3)$ non-empty blocks.  Applying Cauchy's
inequality before the spectral average gives
\begin{equation}\label{eq:all-blocks-final}
 \sum_{\substack{j\ {\rm even}\\T<t_j\le T+U}}\frac1{L_j}
 \left|\sum_{\mathfrak B}P_{j,\mathfrak B}\right|^2
 \ll T^{1+\eps}.
\end{equation}
By Lemmas~\ref{lem:K-tail} and \ref{lem:weighted-truncation},
together with the standard bound $L_j\ll_\eps T^\eps$,
\[
 \sum_{\substack{j\ {\rm even}\\T<t_j\le T+U}}
 L_j|\mathcal E_j^{\mathrm{far}}
  +\mathcal E_j^{\mathrm{shift}}|^2\ll_{A,\psi,\eps}T^{-A},
\]
after choosing the decay exponent in the two tail bounds sufficiently
large and using the spectral counting bound.
Together with \eqref{eq:all-blocks-final}, this proves \eqref{eq:Svariance}.
\end{proof}

\bigskip

\begin{proof}[Proof of Theorem~\ref{thm:main}]
Since
$\mu_j(\psi)-M(\psi)=D_j(\psi)-M(\psi)+S_j(\psi)$,
\[
 |\mu_j(\psi)-M(\psi)|^2
 \le2|D_j(\psi)-M(\psi)|^2+2|S_j(\psi)|^2.
\]

Combining \eqref{eq:Svariance} and \eqref{eq:Dvariance} proves
\[
 \sum_{\substack{j\ {\rm even}\\T<t_j\le T+U}}
 L_j|\mu_j(\psi)-M(\psi)|^2\ll_{\psi,\eps}T^{1+\eps},
\]
and hence Theorem~\ref{thm:main}.
\end{proof}

\appendix
\renewcommand{\theequation}{\thesection.\arabic{equation}}

\section{The Weyl law with prescribed parity}\label{app:parity-weyl}
For completeness, we record the parity-specific Weyl law used in the proof
of Corollary~\ref{cor:rque}.

\begin{lem}\label{prop:parity-weyl}
Let $\{u_j\}$ be an orthonormal basis of Maass cusp forms on
$X=\mathrm{PSL}_2(\mathbb Z)\backslash\mathbb H$, chosen to diagonalise the
reflection $Ru(z)=u(-\overline z)$. Write
\[
 \Delta u_j=(\tfrac14+t_j^2)u_j,\qquad
 Ru_j=\varepsilon_j u_j,\qquad t_j>0,\quad \varepsilon_j\in\{1,-1\}.
\]
Then, for $T\ge2$, counting with multiplicity,
\[
 N_\pm(T):=\#\{j:0<t_j<T,\ \varepsilon_j=\pm1\}
 =\frac{T^2}{24}+O(T\log T).
\]

\end{lem}

\begin{proof}
Put
\[
 N(T)=\sum_{0<t_j<T}1,\qquad
 A(T)=\sum_{0<t_j<T}\varepsilon_j.
\]
The classical Weyl law for the full cuspidal spectrum gives
\begin{equation}\label{eq:full-weyl}
 N(T)=\frac{T^2}{12}+O(T\log T);
\end{equation}
see \cite[Theorem~2.28]{Hejhal}.
Since $N_-(T)=(N(T)-A(T))/2$, it remains to prove
$A(T)\ll T\log T$.

Fix $0<\delta<2\log((1+\sqrt5)/2)$, and choose a nonnegative even
Schwartz function $\rho$ such that
\[
 \int_{\mathbb R}\rho(r)\,dr=1,\qquad
 \supp\widehat\rho\subset(-\delta,\delta),\qquad
 \widehat\rho(u)=\int_{\mathbb R}\rho(r)e^{-iru}\,dr.
\]
 Define
\[
 h_T(r)=\int_{-T}^{T}\rho(r-v)\,dv,\qquad
 g_T(u)=\frac1{2\pi}\int_{\mathbb R}h_T(r)e^{-iru}\,dr
       =\widehat\rho(u)\frac{\sin(Tu)}{\pi u}.
\]
These are admissible test functions for the Selberg trace formula.
Here $\rho$, and hence $h_T$, is understood through its entire
extension furnished by its compactly supported Fourier transform. Indeed, writing $d=\bigl||r|-T\bigr|$, we have
\[
h_T(r)=\int_{r-T}^{r+T}\rho(u)\,du.
\]
If $|r|<T$, the interval $[r-T,r+T]$ contains $[-d,d]$;
if $|r|\ge T$, it is contained in $\{u:|u|\ge d\}$.
Since $\rho\ge0$ and $\int_{\mathbb R}\rho(u)\,du=1$, both cases give
\begin{align}\label{1-approximate}
\bigl|h_T(r)-\mathbf1_{(-T,T)}(r)\bigr|
\le \int_{|u|\ge d}\rho(u)\,du
\ll_B (1+d)^{-B},
\end{align}
where the last estimate follows from the Schwartz decay of $\rho$.

Apply \cite[Theorem~2.1]{BookerLeeStrombergsson2020} with
$N=1$, $\chi=1$ and $n=-1$. The operator $T_{-1}$
is precisely $R$. With the Fourier convention above, the resulting
identity is
\begin{align}\label{eq:reflection-trace}
 h_T(i/2)+\sum_j\varepsilon_jh_T(t_j)
 ={}&\mathcal H_-(g_T)
 -\frac1{2\pi}\int_{\mathbb R}h_T(r)
   \frac{\Gamma'}{\Gamma}(1+ir)\,dr\notag\\
 &+\log(\pi\sqrt2)g_T(0)
 +2\sum_{m\ge2}\frac{\Lambda(m)}m g_T(2\log m).
\end{align}
Here $h_T(i/2)$ is the contribution of the constant eigenfunction, and
substituting \cite[(2.28), (2.34)]{BookerLeeStrombergsson2020} into the
$(\mathrm{NEll}+\mathrm{Ell})$ term of \cite[(2.24)]{BookerLeeStrombergsson2020} gives
\begin{equation}\label{eq:explicit-H}
 \mathcal H_-(g_T)=\sum_{t\in\mathbb{Z}\setminus\{0\}}c_tg_T(\ell_t),
\end{equation}
where
\[
 \ell_t=2\log\frac{|t|+\sqrt{t^2+4}}2,\qquad
 t^2+4=d_tf_t^2,
\]
$d_t$ is a fundamental discriminant, $f_t\ge1$, and
\[
 c_t=\frac{h^+(d_t)\log\eta_{d_t}}{\sqrt{t^2+4}}
 \prod_{p^a\parallel f_t}
 \left(1+\left(p-\left(\frac{d_t}{p}\right)\right)
                  \frac{p^a-1}{p-1}\right).
\]
In this expression, $h^+(d)$ denotes the narrow ideal class number of
$\mathbb Q(\sqrt d)$, $\eta_d>1$ is the least unit greater than $1$
of norm $1$, and $(\frac d p)$ is the Kronecker symbol.
For every $t\neq 0$,
\[
\ell_t\ge2\log\frac{1+\sqrt5}2>\delta.
\]
Since $\supp g_T\subset(-\delta,\delta)$,
\eqref{eq:explicit-H} gives $\mathcal H_-(g_T)=0$.
The prime-power sum in \eqref{eq:reflection-trace} also vanishes, since
$2\log m\ge2\log2>\delta$. Stirling's formula, together with
$g_T(0)=T/\pi$ and $h_T(i/2)=O(1)$, now yields
\begin{equation}\label{eq:smoothed-parity}
 \sum_j\varepsilon_jh_T(t_j)\ll
 \int_{\mathbb R}h_T(r)\log(2+|r|)\,dr+T+1
 \ll T\log T.
\end{equation}

To remove the smoothing, \eqref{eq:full-weyl} implies
\[
 \#\{j:k\le t_j<k+1\}\ll(k+1)\log(k+2)
 \qquad(k\in\mathbb{Z}_{\ge0}).
\]
The estimate \eqref{1-approximate} with $B>3$ then gives
\[
 \sum_j\bigl|h_T(t_j)-\mathbf1_{(0,T)}(t_j)\bigr|
 \ll T\log T.
\]
Consequently \eqref{eq:smoothed-parity} implies $A(T)\ll T\log T$,
which proves the lemma.
\end{proof}


\begin{thebibliography}{99}






\bibitem{BookerLeeStrombergsson2020}
A.~R. Booker, M. Lee and A. Str\"ombergsson,
Twist-minimal trace formulas and the Selberg eigenvalue conjecture,
J. Lond. Math. Soc. (2) {\bf 102} (2020), no.~3, 1067--1134.

\bibitem{FaaDiBruno1857}

F.~Fa{\`a} di Bruno, Note sur une nouvelle formule de calcul
diff{\'e}rentiel, \emph{Quart. J. Pure Appl. Math.} \textbf{1} (1857),
359--360.

\bibitem{GhoshReznikovSarnak2013}
A. Ghosh, A. Reznikov and P.~C. Sarnak, Nodal domains
of Maass forms I, Geom. Funct. Anal. {\bf 23} (2013),
no.~5, 1515--1568.



\bibitem{GradshteynRyzhik2007}
I.~S. Gradshteyn and I.~M. Ryzhik, \emph{Table of integrals,
series, and products}, 7th ed., translated from the Russian; translation
edited and with a preface by A. Jeffrey and D. Zwillinger,
Elsevier/Academic Press, Amsterdam, 2007.

\bibitem{Hejhal}
D.~A. Hejhal, {\it The Selberg trace formula
for ${\rm PSL}(2,\,{\bf R})$. Vol. 2},
Lecture Notes in Mathematics, 1001, Springer, Berlin, 1983.

\bibitem{HoffsteinLockhart1994}
J.~E. Hoffstein and P. Lockhart, Coefficients
of Maass forms and the Siegel zero,
Ann. of Math. (2) {\bf 140} (1994), no.~1, 161--181.

\bibitem{Iwaniec2002}
H. Iwaniec, \emph{Spectral methods of automorphic forms}, 2nd ed.,
Graduate Studies in Mathematics, 53, Amer. Math. Soc., Providence, RI;
Revista Matem\'atica Iberoamericana, Madrid, 2002.

\bibitem{Jutila2000}
M.~I. Jutila, On spectral large sieve inequalities,
Funct. Approx. Comment. Math. {\bf 28} (2000), 7--18.

\bibitem{Li2010}
X. Li, Upper bounds on $L$-functions at the edge of the
critical strip, Int. Math. Res. Not. IMRN {\bf 2010},
no.~4, 727--755.


\bibitem{Lindenstrauss2006}
E. Lindenstrauss, Invariant measures and
arithmetic quantum unique ergodicity,
Ann. of Math. (2) {\bf 163} (2006), no.~1, 165--219.

\bibitem{Liu2017}
S. Liu, {\it Automorphic L-Functions and Their Derivatives},
ProQuest LLC, Ann Arbor, MI, 2017.

\bibitem{LuoSarnak1995}
W. Luo and P.~C. Sarnak, Quantum ergodicity
of eigenfunctions on ${\rm PSL}_2(\mathbb Z)\backslash \mathbb H^2$,
Inst. Hautes \'Etudes Sci. Publ. Math. No. 81 (1995), 207--237.



\bibitem{Soundararajan2010}
K. Soundararajan, Quantum unique ergodicity
for ${\rm SL}_2(\mathbb Z)\backslash\mathbb H$,
Ann. of Math. (2) {\bf 172} (2010), no.~2, 1529--1538.

\bibitem{SunZhang2026}
Q. Sun and Q. Zhang, Mass distribution for holomorphic
cusp forms on the vertical geodesic,
J. Number Theory {\bf 279} (2026), 827--857.

\bibitem{Teschl2012}
G. Teschl, {\it Ordinary differential equations
and dynamical systems}, Graduate Studies in Mathematics, 140,
Amer. Math. Soc., Providence, RI, 2012.

\bibitem{Young2016}
M.~P. Young, The quantum unique ergodicity conjecture
for thin sets, Adv. Math. {\bf 286} (2016), 958--1016.

\bibitem{Zenz2021}
P. Zenz, Quantum variance for holomorphic Hecke cusp forms on the
vertical geodesic, arXiv:2111.04713, 2021.

\end{thebibliography}
\end{document}